\documentclass[11pt]{amsart}
\usepackage{amsmath}
\usepackage{epic,eepic}			

\newcommand{\bi}{\bf}                        

\numberwithin{equation}{section}   

\newcommand{\la}{\langle}               
\newcommand{\ra}{\rangle}               
\renewcommand{\(}{\left(}
\renewcommand{\)}{\right)}

\renewcommand{\]}{\right]}

\renewcommand{\Re}{{\mathbf R}}             
\newcommand{\e}{\varepsilon}             

\newcommand{\nn}{\nonumber}             
\newcommand{\note}[1]{}                    
\newcommand{\commentout}[1]{}

\newtheorem{thm}{Theorem}[section]
\newtheorem{lemma}[thm]{Lemma}
\newtheorem{prop}[thm]{Proposition}
\newtheorem{cor}[thm]{Corollary}

\theoremstyle{definition}
\newtheorem{definition}[thm]{Definition}

\theoremstyle{remark}
\newtheorem{remark}[thm]{Remark}

\newcommand{\f}{\partial}

\newcommand{\nor}{{\mathbf n}}

\renewcommand{\phi}{\varphi}

\newcommand{\trace}{\operatorname{trace}}
\newcommand{\cd}{,\dots,}
\newcommand{\ric}{\operatorname{Ric}\nolimits}
\newcommand{\Ric}{\operatorname{Ric}\nolimits}
\newcommand{\mean}{{\mathcal M}}

\newcommand{\ecst}{C_{\mathbf E}}
\newcommand{\scst}{C_{\mathbf S}}
\newcommand{\hcst}{C_{\mathbf H}}
\newcommand{\ol}{\overline}

\newcommand{\els}{{\mathbf S}}
\newcommand{\hyps}{{\mathbf H}}

\title[Maximum Principle and Lorentzian Geometry]{A Strong
 Maximum Principle for Weak Solutions of Quasi-Linear Elliptic
Equations with Applications to Lorentzian and Riemannian Geometry}
\author[Andersson]{Lars Andersson\footnotemark{$^*$}}
\thanks{\footnotemark{$^*$} Supported in part by NFR, contract no. F-FU
4873-307.}
\address{Department of Mathematics \\
Royal Institute of Technology \\
S-100 44 Stockholm, Sweden}
\email{larsa\char'100math.kth.se}

\author[Galloway]{Gregory J. Galloway$^\dagger$}
\thanks{\footnotemark{$^\dagger$} Supported in part by NSF grant
DMS-9204372}
\address{Department of Mathematics \\
University of Miami\\
Coral Gables, FL 33124, USA}
\email{galloway\char'100math.miami.edu }

\author[Howard]{Ralph Howard$^\ddagger$}
\thanks{\footnotemark{$^\ddagger$} Supported in part by DEPSCoR grants
N00014-94-1-1163 and DAAH-04-96-1-0326}
\address{Department of Mathematics\\
University of South Carolina \\
Columbia, S.C. 29208, USA}
\email{howard\char'100math.sc.edu}


\begin{document}

\begin{abstract}
The strong maximum principle is proved to hold 
for weak (in the sense of support functions) sub- and super-solutions to a 
class of quasi-linear elliptic equations that includes the mean curvature 
equation for
$C^0$~spacelike hypersurfaces in a Lorentzian manifold.  As one
application a Lorentzian warped product splitting theorem is given.
\end{abstract}
\maketitle
\tableofcontents


\section{Introduction}\label{sec:intro}

The geometric maximum principle for $C^2$ hypersurfaces, first proven
by Alexandrov~\cite{Alexandrov:uniqueness1,Alexandrov:uniqueness2},
has become a basic tool in differential geometry especially for
proving uniqueness and rigidity results.  The geometric maximum
principle is a consequence of the strong maximum principle of
Alexandrov~\cite{Alexandrov:2nd-thms} for $C^2$ sub- and
super-solutions for nonlinear uniformly elliptic equations which is in
turn based on the original maximum principle of
E.~Hopf~\cite{Hopf:maximum} for linear elliptic equations.

However there are  naturally occurring geometric situations, for
example level sets of Busemann functions, where the
classical maximum principle does not apply as the hypersurfaces
involved are only of class $C^0$ (that is locally graphs of continuous
functions) and only satisfy curvature inequalities in the sense of
support functions (see Definitions~\ref{mean-spt},
\ref{def:hyper-mean}, and \ref{mean-hyper-def}).  In this paper we
give a version of the strong maximum principle, Theorem~\ref{thm:max},
and deduce geometric maximum principles, Theorems~\ref{geo:max}
and~\ref{Riem-max}, general enough to cover most recent applications of
maximum principles to uniqueness and rigidity questions in Riemannian
and Lorentzian geometry.

Our main analytic result is Theorem~\ref{thm:max} which extends the
strong maximum principle to weak (in the sense of support functions)
sub- and super-solutions of uniformly elliptic quasi-linear equations.
For linear equations this is due to Calabi~\cite{Calabi:maximum} who
also introduced the idea of sub-solutions in the sense of support
functions.  The geometric versions are a strong maximum principle,
Theorem~\ref{geo:max}, for $C^0$~spacelike hypersurfaces
(Definition~\ref{def:spacelike}) in Lorentzian manifolds and
Theorem~\ref{Riem-max} which is a maximum principle for
hypersurfaces in Riemannian manifolds that can be locally represented
as graphs. Eschenburg~\cite{Eschenburg:max} gives a version of these
geometric results under the extra assumption that one of the two
hypersurfaces is $C^2$, but in the applications given here neither
hypersurface will have any {\em a priori\/} smoothness.

A very natural example of rough hypersurfaces where our maximum
principles work well is the level sets of Busemann functions in
Riemannian or Lorentzian manifolds.  By applying our results to these
level sets in Lorentzian manifolds we prove a warped product splitting
theorem which can be viewed as an extension of the Lorentzian
splitting
theorem~\cite{Eschenburg:split,galloway:split,newman:splitting} to
spacetimes that satisfy the strong energy condition with a positive
cosmological constant.  (See
\S\ref{sec:appl:warp} 
for a discussion of how proofs of splitting theorems and warped
product splitting theorems are simplified by use of the maximum
principle for rough hypersurfaces.)  As applications of warped product
splitting we obtain a characterization, Corollary~\ref{Lor-max-diam},
of the universal anti-de~Sitter space which can be viewed as a
Lorentzian analogue of the Cheng maximum diameter
theorem~\cite{Cheng:max-diam}.  Certain
Lorentzian warped products that are locally spatially isotropic are
characterized (Theorem~\ref{conformal}) by the existence of a line
together with a version of the Weyl curvature hypothesis of
R.~Penrose. Further applications of  maximum principles for rough
hypersurfaces are given in
\cite{andersson:howard:RW,galloway:warsaw,Howard:inner-sphere}.

\section{The Analytic Maximum Principle}
\label{sec:analytmax}

We  denote partial derivatives on $\Re^n$ (with coordinates $x^1\cd
x^n$) by
$$
D_i f:=\frac{\f f}{\f x^i},\qquad D_{ij}f:=\frac{\f^2 f}{\f x^i\f x^j} .
$$
Also let
$$
Df:=(D_1f\cd D_nf),\qquad D^2f:=[D_{ij}f].
$$
Thus $Df$ is the gradient of $f$ and $D^2f$ is the matrix of second order
partial derivatives of $f$, i.e. the Hessian of $f$.  For
each point $x$ in the domain of $f$ the Hessian $D^2f(x)$ is a
symmetric matrix.  
If $A$ and $B$ are symmetric matrices
we write $A\le B$ if $B-A$ is positive semi-definite.

As coordinates on $\Re^n\times\Re\times\Re^n$ use $(x,r,p)=(x^1\cd x^n,r,p^1\cd
p^n)$.  If $U\subseteq\Re^n\times\Re\times\Re^n$ the
{\bi fiber\/} over $x\in \Re^n$ is
$$
U_x:=\{(r,p)\in\Re\times\Re^n: (x,r,p)\in U\}.
$$
Note that the fiber may be empty.  
\begin{definition}\label{def-addmissible}
Let $\Omega$ be an open set and
$U\subset\Re^n\times\Re\times\Re^n$.  Then a $C^1$ function $u:\Omega\to\Re$ is
{\bi $U$--admissible\/} over $\Omega$ if and only if for all $x\in\Omega$,
$(x,u(x),D u(x))=(x,u(x),D_1u(x)\cd D_nu(x))\in U$.  (In particular this
implies the fiber $U_x\ne\emptyset$ for all $x\in \Omega$.)
~\qed
\end{definition}
Let
$U\subseteq\Re^n\times\Re\times\Re^n$ be open.  Then a {\bi quasi-linear
operator\/} $\mean$ on $U$ is a collection of functions $a^{ij}=a^{ji},
b:U\to\Re$ $1\leq i,j \leq n$.  
If $u$ is a $U$--admissible $C^2$ function defined on the open
set $\Omega\subseteq\Re^n,$ then $\mean$ is evaluated on $u$ by
$$
\mean[u]=\sum_{i,j}a^{ij}(x,u, D u) D_{ij}u+b(x,u, D u).
$$

\begin{definition}
\label{uniform}
The quasi-linear operator $\mean$ on the open set
$U\subseteq \Re^n\times\Re\times\Re^n$ is {\bi uniformly elliptic\/} on $U$
if and only if
\begin{enumerate}
\item The functions $a^{ij}$ and $b$ are continuous on $U$ and are $C^1$
functions of the arguments~$(r,p)$.
\item There is a constant $\ecst$ so that for all $(x,r,p)\in U$
\begin{equation}\label{a-est}
\frac1{\ecst}\|\xi\|^2\le \sum_{i,j}a^{ij}(x,r,p)\xi_i\xi_j\le
\ecst\|\xi\|^2
\end{equation}
and
\begin{equation}\label{der-bd}
\left|\frac{\f a^{ij}}{\f p^k}\right|,\ \left|\frac{\f a^{ij}}{\f
r}\right|,\
\left|\frac{\f b}{\f p^k}\right|,\  \left|\frac{\f b}{\f r}\right|,\ |b|\le
\ecst.
\end{equation}
The constant $\ecst$ is the {\bi constant of ellipticity\/} of $\mean$ on
$U$.
\end{enumerate}
~\qed
\end{definition}

We will need the notion  of a support function.  Let $u$ be defined on the
open set $\Omega\subset\Re^n$ and let $x_0\in \Omega$.  Then $\phi$ is an
{\bi upper support function\/} for $u$ at $x_0$ if and only if
$\phi(x_0)=u(x_0)$ and for some open neighborhood $N$ of $x_0$ the
inequality~$\phi\ge u$ holds.  The function $\phi$ is a {\bi lower support
function\/} for $u$ at $x_0$ if and only if $\phi(x_0)=u(x_0)$ and $\phi\le
u$ in some open neighborhood $N$ of $x_0$.
Given a sequence $B^1\cd B^n$ of functions  set
$$
\|B\|(x):=\sqrt{ (B^1(x))^2+\cdots+(B^n(x))^2}
$$
and if $\phi$ is a $C^2$ function we set
$$
|D^2\phi|(x):=\max_{i,j}|D_{ij}\phi(x)|.
$$

\begin{definition}\label{mean-spt}
Let $U\subseteq\Re^n\times\Re\times\Re^n$ be open, let $\mean$ be a
uniformly 
elliptic operator on $U$ and let $H_0$ be a constant.
\begin{enumerate}\item
A lower semi-continuous function $u$
defined on the open subset $\Omega\subseteq\Re^n$  satisfies
$\mean[u]\le H_0$ {\bi in the sense of support functions}\/ iff
for all $\e>0$ and  $x\in \Omega$ there is a $U$--admissible $C^2$ upper 
support function $\phi_{x,\e}$ to $u$ at $x$ 
so that
$$
\mean[\phi_{x,\e}](x) \le H_0+\e.
$$
\item
An upper semi-continuous function $v$ defined on $\Omega$
satisfies $\mean[v]\ge H_0$ in the {\bi sense of support
functions with a one--sided Hessian bound}\/ iff  there is a constant
$\scst >0$ (the support
constant) so that for all
$\e>0$ and $x\in\Omega$ there is a $U$--admissible $C^2$ lower support function
$\psi_{x,\e}$ for $v$ at $x$ so that
$$
\mean[\psi_{x,\e}](x)\ge H_0-\e,\quad\text{and}\quad
D^2\psi_{x,\e}(x)\ge-\scst I
$$
(The constant $\scst$ can depend on $v$ but is independent of $\e$ and $x$.)
\end{enumerate}
~\qed
\end{definition}

We are now ready to state the main theorem of this section.
\begin{thm}[Maximum Principle]\sl\label{thm:max}
Let $\Omega\subseteq\Re^n$ be open and connected and let $U\subseteq
\Re^n\times\Re\times\Re^n$ be such that for each $x\in \Omega$, the
fiber $U_x$ is nonempty and convex.  Let $u_0,u_1:\Omega\to\Re$ be
functions satisfying the following. 

\begin{enumerate}
\item $u_0$ is lower semi-continuous and $\mean[u_0]\le H_0$ in the sense
of support functions.
\item $u_1$ is upper-semi-continuous and $\mean[u_1]\ge H_0$ in the sense
of support functions with a one--sided Hessian bound.
\item $u_1\le u_0$ in $\Omega$, and there is at least one point
$x\in\Omega$ with $u_1(x)=u_0(x)$.
\end{enumerate}
Then $u_0\equiv u_1$ in $\Omega$ and $u_0=u_1$ is locally a $C^{1,1}$ function
in $\Omega$.  Finally if $a^{ij}$ and $b$ are locally $C^{k,\alpha}$ functions
in $U$ then $u_0=u_1$ is locally a $C^{k+2,\alpha}$ function in
$\Omega$.  In
particular if $a^{ij}$ and $b$ are smooth, then so is $u_0=u_1$.
\end{thm}

\begin{remark}
The geometric  maximum principle for $C^2$ sub- and super-solutions to
the Laplacian $\Delta:=\sum_{i}D_{ii}$ is classical and follows from
the fact that the difference $f=u_1-u_0$ is sub-harmonic and use of the
sub-mean value theorem for sub-harmonic functions.  For
$C^2$ sub- and super-solutions to general variable coefficient linear
operators the result is due to E.~Hopf~\cite{Hopf:maximum}.  
The extension to nonlinear uniformly elliptic equations was done by
Alexandrov~\cite{Alexandrov:2nd-thms,Alexandrov:uniqueness2}.
The formulation and proof of the strong maximum principle for linear
operators with sub- and super-solutions  defined in the sense of
support functions is due to Calabi~\cite{Calabi:maximum}.  
A partial geometric extension of Calabi's result to the mean 
curvature operator of
hypersurfaces in Riemannian and Lorentzian manifolds was given by
Eschenburg~\cite{Eschenburg:max} who also introduced  (in a geometric
context) the idea of $\mean[u_0]\ge H_0$ in the sense of support
functions with a one sided Hessian bound.
%
%
~\qed
\end{remark}
The proof of Theorem \ref{thm:max} will be given in \S \ref{proof}. 
An outline of the proof is given in \S \ref{sec:outline} and some 
preliminaries are given in \S \ref{setup} and \S \ref{sec:lemmata}.

\subsection{Outline of the Proof}\label{sec:outline}
Here we give the ideas behind the proof which, stripped of the
technical details, are easy.  The basic outline of the proof
and the idea of using support functions is as in
Calabi~\cite{Calabi:maximum} which in turn owes much to the original paper
of E.~Hopf~\cite{Hopf:maximum}.   

First let $\mean$ be as in
Theorem~\ref{thm:max} and let $\phi_0$ and $\phi_1$ be $C^2$ functions.
Then (Lemma~\ref{lemma:coeff} below) there are functions
$A^{ij}[\phi_0,\phi_1]$, $B^i[\phi_0,\phi_1]$, $C[\phi_0,\phi_1]$ so that
\begin{align}
\mean[\phi_1]-\mean[\phi_0]
&=\sum_{i,j}A^{ij}[\phi_0,\phi_1]D_{ij}(\phi_1-\phi_0)
        +\sum_iB^i[\phi_0,\phi_1]D_i(\phi_1-\phi_0)\nn\\
\label{com}     &\qquad +C[\phi_0,\phi_1](\phi_1-\phi_0)
\end{align}
and such that the matrix $\[A^{ij}[\phi_0,\phi_1]\]$ is pointwise
positive definite.
Given $\phi_0$ and $\phi_1$ define the linear differential
operator
$L_{\phi_0,\phi_1}$
\begin{equation}\label{Lbar}
{L}_{\phi_0,\phi_1}:=\sum_{i,j}A^{ij}[\phi_0,\phi_1]D_{ij}
        +\sum_i{B}^i[\phi_0,\phi_1]D_i.
\end{equation}

With this notation and the hypotheses of Theorem~\ref{thm:max}  assume,
toward a contradiction, that $u_0\ne u_1$.
Then $K:=\{x\in \Omega :u_1(x)=u_0(x)\}$ is a proper closed subset of
$\Omega$.  
\begin{figure}
\centering

\setlength{\unitlength}{0.00050000in}
\begingroup\makeatletter\ifx\SetFigFont\undefined
\def\x#1#2#3#4#5#6#7\relax{\def\x{#1#2#3#4#5#6}}%
\expandafter\x\fmtname xxxxxx\relax \def\y{splain}%
\ifx\x\y   
\gdef\SetFigFont#1#2#3{%
  \ifnum #1<17\tiny\else \ifnum #1<20\small\else
  \ifnum #1<24\normalsize\else \ifnum #1<29\large\else
  \ifnum #1<34\Large\else \ifnum #1<41\LARGE\else
     \huge\fi\fi\fi\fi\fi\fi
  \csname #3\endcsname}%
\else
\gdef\SetFigFont#1#2#3{\begingroup
  \count@#1\relax \ifnum 25<\count@\count@25\fi
  \def\x{\endgroup\@setsize\SetFigFont{#2pt}}%
  \expandafter\x
    \csname \romannumeral\the\count@ pt\expandafter\endcsname
    \csname @\romannumeral\the\count@ pt\endcsname
  \csname #3\endcsname}%
\fi
\fi\endgroup
{\renewcommand{\dashlinestretch}{30}
\begin{picture}(5942,4174)(0,-10)
\put(2775,1807){\ellipse{3600}{3600}}
\put(2775,1807){\ellipse{2400}{2400}}
\put(2775,3007){\ellipse{1050}{1050}}
\path(2775,3007)(3375,2407)
\path(2881.066,2943.360)(2775.000,3007.000)(2838.640,2900.934)
\path(3375,2407)(5175,2407)
\path(3300,3082)(5175,3082)
\path(3420.000,3112.000)(3300.000,3082.000)(3420.000,3052.000)
\path(2400,2632)(1875,2182)
\path(2328.413,2531.127)(2400.000,2632.000)(2289.365,2576.683)
\path(1875,2182)(375,2182)
\path(2775,1807)(5175,1807)
\path(2895.000,1837.000)(2775.000,1807.000)(2895.000,1777.000)
\path(2625,2857)(375,2857)
\path(2505.000,2827.000)(2625.000,2857.000)(2505.000,2887.000)
\path(1125,3982)	(1191.234,3988.486)
	(1252.612,3993.113)
	(1361.676,3996.573)
	(1453.950,3991.940)
	(1531.193,3978.775)
	(1647.613,3925.089)
	(1725.000,3832.000)

\path(1725,3832)	(1686.769,3713.571)
	(1650.000,3607.000)

\path(1650,3607)	(1724.860,3572.470)
	(1827.986,3588.291)
	(1935.870,3613.467)
	(2025.000,3607.000)

\path(2025,3607)	(2089.209,3547.470)
	(2143.057,3459.209)
	(2194.127,3369.843)
	(2250.000,3307.000)

\path(2250,3307)	(2315.702,3281.863)
	(2400.068,3267.010)
	(2484.399,3253.402)
	(2550.000,3232.000)

\path(2550,3232)	(2637.806,3160.934)
	(2697.593,3096.962)
	(2775.000,3007.000)

\path(2775,3007)	(2885.066,3067.987)
	(2967.068,3110.972)
	(3075.000,3157.000)

\path(3075,3157)	(3190.117,3175.853)
	(3261.346,3181.046)
	(3336.547,3185.661)
	(3411.868,3191.381)
	(3483.454,3199.889)
	(3600.000,3232.000)

\path(3600,3232)	(3710.096,3345.846)
	(3825.000,3457.000)

\path(3825,3457)	(3897.084,3469.777)
	(3983.570,3461.624)
	(4079.014,3441.180)
	(4177.969,3417.085)
	(4274.990,3397.980)
	(4364.633,3392.504)
	(4441.451,3409.297)
	(4500.000,3457.000)

\path(4500,3457)	(4508.681,3541.026)
	(4479.265,3601.758)
	(4425.000,3682.000)

\path(4425,3682)	(4450.377,3753.486)
	(4468.504,3815.742)
	(4483.005,3916.964)
	(4468.504,3994.453)
	(4425.000,4057.000)

\path(4425,4057)	(4349.018,4110.556)
	(4262.873,4136.210)
	(4170.449,4140.261)
	(4075.628,4129.011)
	(3982.290,4108.760)
	(3894.318,4085.807)
	(3815.594,4066.454)
	(3750.000,4057.000)

\path(3750,4057)	(3660.870,4050.613)
	(3552.246,4039.251)
	(3430.536,4024.976)
	(3366.777,4017.390)
	(3302.149,4009.847)
	(3237.454,4002.607)
	(3173.493,3995.925)
	(3050.977,3985.270)
	(2941.010,3979.941)
	(2850.000,3982.000)

\path(2850,3982)	(2783.585,3996.949)
	(2702.261,4023.704)
	(2619.807,4048.357)
	(2550.000,4057.000)

\path(2550,4057)	(2462.564,4031.605)
	(2362.500,3982.000)
	(2262.436,3932.395)
	(2175.000,3907.000)

\path(2175,3907)	(2104.477,3913.406)
	(2020.830,3934.330)
	(1939.268,3960.339)
	(1875.000,3982.000)

\path(1875,3982)	(1792.967,4015.512)
	(1691.017,4062.512)
	(1587.309,4106.757)
	(1500.000,4132.000)

\path(1500,4132)	(1415.703,4138.669)
	(1306.542,4139.741)
	(1183.598,4135.721)
	(1120.421,4131.961)
	(1057.954,4127.118)
	(940.690,4114.437)
	(842.889,4098.185)
	(750.000,4057.000)

\path(750,4057)	(831.176,4023.104)
	(947.113,4004.332)
	(1027.763,3993.721)
	(1125.000,3982.000)

\put(5250,1732){\makebox(0,0)[lb]{\smash{{{\SetFigFont{11}{13.2}{rm}$x_0$}}}}}
\put(5250,2332){\makebox(0,0)[lb]{\smash{{{\SetFigFont{11}{13.2}{rm}$x_1$}}}}}
\put(5250,3007){\makebox(0,0)[lb]{\smash{{{\SetFigFont{11}{13.2}{rm}$S''$}}}}}
\put(75,2107){\makebox(0,0)[lb]{\smash{{{\SetFigFont{11}{13.2}{rm}$S'$}}}}}
\put(3150,3682){\makebox(0,0)[lb]{\smash{{{\SetFigFont{11}{13.2}{rm}$K$}}}}}
\put(0,2857){\makebox(0,0)[lb]{\smash{{{\SetFigFont{11}{13.2}{rm}$x_*$}}}}}
\put(2100,307){\makebox(0,0)[lb]{\smash{{{\SetFigFont{11}{13.2}{rm}$B(x_0,3r_0)$}}}}}
\put(2100,982){\makebox(0,0)[lb]{\smash{{{\SetFigFont{11}{13.2}{rm}$B(x_0,2r_0)$}}}}}
\end{picture}
}

\caption[]{}
\label{fig1}
\end{figure}
Thus there is a closed ball $\ol{B}(x_0,2r_0)$ so that
$\ol{B}(x_0,3r_0)\subset\Omega$ and $\ol{B}(x_0,2r_0)$ meets $K$ in exactly
one point $x_1$ and this point is on $\f B(x_0,2r_0)$.  Let $r_1\le r_0$
and let
$$
S':=\f B(x_1,r_1)\cap\ol{B}(x_0,2r_0), \quad
S'':=\f B(x_1,r_1)\setminus
\ol{B}(x_0,2r_0),
$$
see Figure \ref{fig1}.
Then $(u_1-u_0)<0$ on $S'$ as $S'$ is disjoint from $K$.  For
$\alpha>0$ set $w(x):=\|x-x_0\|^{-\alpha}-\|x_1-x_0\|^{-\alpha}
=\|x-x_0\|^{-\alpha}-(2r_0)^{-\alpha}$.  As $w$ is bounded on the
compact set $S'$ and $(u_1-u_0)$ is negative on $S'$ it is possible to
choose $\delta>0$ so that if $h$ is the function $h:= (u_1-u_0)+\delta
w$ then $h<0$ on $S'$.  But $w$ is a decreasing function of
$\|x-x_0\|$ (so $w(x)< 0$ on $S''$) and $(u_1-u_0)\le 0$ so also $h<0$
on $S''$.
Thus $h<0$ on $S'\cup S''=\f B(x_1,r_1)$ and
$h(x_1)=(u_1(x_1)-u_0(x_0))+\delta w(x_1)=0$, so the function
$h$ has a nonnegative interior maximum in $\ol{B}(x_1,r_1)$ at a point $x_*$.
Let $\phi_0$ be an upper support function for $u_0$ at $x_*$ and
$\phi_1$ a lower support function for $u_1$ at $x_*$ so that for some
very small $\e$ the inequalities $\mean[\phi_1]\ge H_0-\e$ and
$\mean[\phi_0]\le H_0+\e$ hold.  Then~(\ref{com}) and~(\ref{Lbar})
imply
\begin{align*}
L_{\phi_0,\phi_1}(\phi_1-\phi_0)&\ge
                -2\e-|C[\phi_0,\phi_1]||\phi_1(x_*) -\phi_0(x_*)|\\
   &=-2\e-|C[\phi_0,\phi_1]||u_1(x_*)-u_0(x_*)|.
\end{align*}
Choose $\alpha$ in
the definition of $w$ so that $L_{\phi_0,\phi_1}w>0$.
The inequalities $h(x_*) \ge 0$ and $u_1(x_*)-u_0(x_*) \le 0$ imply that
$|u_1(x_*)-u_0(x_*)| \le \delta |w(x_*)|$.  As $w(x_1) = 0$,
by choosing $r_1$ small enough we make the inequality
$$
|C[\phi_0,\phi_1]||u_1(x_*)-u_0(x_*)|<\e .
$$
hold. 
Then $f:=(\phi_1-\phi_0)+\delta w$ will satisfy
\begin{equation}\label{eq:good}
\begin{split}
L_{\phi_0,\phi_1}f(x_*)&=L_{\phi_0,\phi_1}(\phi_1-\phi_0)(x_*)+\delta
                L_{\phi_0,\phi_1}w(x_*)\\
        &\ge-2\e-|C[\phi_0,\phi_1]||u_1(x_*)-u_0(x_*)|
                +\delta L_{\phi_0,\phi_1}w(x_*)\\
        &\ge-3\e +\delta L_{\phi_0,\phi_1}w(x_*) >0
\end{split}
\end{equation}
provided $\e$ was chosen small enough.  But as $\phi_0$ is an upper
support function for $u_0$ at $x_*$ and $\phi_1$ is a lower support
function for $u_1$ at $x_*$ the function $f=(\phi_1-\phi_0)+\delta w$
will have a local maximum at $x_*$ because $h=(u_1-u_0)+\delta w$ has
a maximum at $x_*$.  But then $Df(x_*)=0$ and $D^2f(x_*)$ is negative
semi-definite and $\[A^{ij}[\phi_0,\phi_1]\]$ is positive definite
which implies $L_{\phi_0,\phi_1}f(x_*)\le0$.  
We have arrived at a contradiction to the inequality (\ref{eq:good})
and this completes the outline of the proof of Theorem \ref{thm:max}.
\bigskip

The proof is more complicated than this outline indicates for two
reasons.  First to show $L_{\phi_0,\phi_1}w(x_*)>0$ and
$ |C[\phi_0,\phi_1]||u_1(x_*)-u_0(x_*)|<\e$ we need bounds on
$B^i[\phi_0,\phi_1]$, and $C[\phi_0,\phi_1]$.  But the definitions of
$B^i[\phi_0,\phi_1]$, and $C[\phi_0,\phi_1]$ involve  the Hessians
$D^2\phi_0$ and $D^2\phi_1$ and therefore we need to find
{\it a~priori}\/ bounds for
$|D^2\phi_0|(x_*)$ and $|D^2\phi_1|(x_*)$.  This is done in
Lemma~\ref{Hess-bd} and its Corollary~\ref{Hess-bd2}
by using that at $x_*$ the function
$f=(\phi_1-\phi_0)+\delta w$ has a maximum at $x_*$ so the
Hessian $D^2f(x_*)$ will be negative semi-definite.  Using
$D^2f(x_*)\le 0$ and the one--sided bound
$D^2\phi_1\ge -\scst I$ gives a lower bound $D^2\phi_0\ge D^2\phi_1+\delta
D^2w$ for $D^2\phi_0$ at $x_*$.  The operator
$\mean$ is uniformly elliptic so the lower bound on $D^2\phi(x_*)$ and
$\mean[\phi_0]\le H_0+\e$ implies an upper bound on $D^2\phi_0$.  Then again
using that $D^2f=(\phi_1-\phi_0)+\delta w\le 0$ so that $D^2\phi_1\le
D^2\phi_0 -\delta D^2w$ we get an upper bound on $D^2\phi_1$.  This gives
the required (two--sided) bounds on
$D^2\phi_0$ and $D^2\phi_1$ in terms of $\delta$ and
$\alpha$ (which appears in the definition of $w$).  Which brings us to
the second problem.  In the outline above the choices of $\alpha$,
$\delta$, $r_1$ etc.\ depend on the bounds of
$B^i[\phi_0,\phi_1]$, $C[\phi_0,\phi_1]$
which in turn depend on $\alpha$, $\delta$, $r_1$ etc.  Thus at
least formally the argument is circular.  However when the relations
required between the various parameters are written out explicitly it
is possible, with some care, to make the choices so that the argument
works.

\subsection{Reduction to a Standard Setup} \label{setup}
We first note if we can prove the theorem under the assumption that
$\Omega$ is convex then we can express an arbitrary connected domain
as a union of convex subdomains and deduce $u_0\equiv u_1$ in the
entire domain by doing a standard ``analytic continuation'' argument.
Thus assume $\Omega$ is convex.  While in some applications it is convenient
to work with the constant $H_0$ in Theorem~\ref{thm:max} being non-zero, by
replacing $b(x,r,p)$ by $b(x,r,p)-H_0$ it is enough to prove the result
with $H_0=0$. (As we are assuming as part of the definition of
uniform ellipticity that $|b(x,r,p)|\le \ecst$ we may have to replace the
constant $\ecst$ by $\ecst +|H_0|$.)  As we will be using the same
hypothesis in several of the lemmata below it is convenient to give a
name to data $x_0$, $x_1$, $r_0$, $r_1$, $\alpha$,
$\delta$, $\phi_0$, $\phi_1$ that satisfies these hypotheses and at
the same time make some normalizations that will simplify notation.
We will choose $x_0$ as above so that the closed ball
$\ol{B}(x_0,2r_0)$ meets $K:=\{ x\in \Omega: u_0(x)=u_1(x)\}$ at
exactly one point $x_1$ and this point is on $\f B(x_0,2r_0)$.

By translating the coordinates we can assume $x_0=0$, the origin of
the coordinate system. With these simplifications we make the following
definition.
\begin{definition}
\label{stand} 
The data
$x_0$, $x_1$, $x_*$, $r_0$, $r_1$, $\alpha$, $\delta$,
$\phi_0$, $\phi_1$ will be said to satisfy the {\bi standard setup}\/
if and only if
\begin{enumerate}
\item $x_0=0$ and $x_1$ is a point with
$\|x_1-x_0\|=\|x_1\|=2r_0$.
\item The positive numbers $r_0$, $r_1$ satisfy
$$
r_1\le r_0<3r_0\le 1.
$$
\item The number $\alpha$ is positive and we use it to define the
{\bi comparison function}\/
$$
w(x):=\|x\|^{-\alpha}.
$$
\item The number $\delta$ is  positive, the point $x_*$ is
in the interior of $B(x_1,r_1)$, the functions $\phi_0$ and $\phi_1$
are $C^2$ and defined in some neighborhood of $x_*$ (but not
necessarily in all of $B(x_1,r_1)$) and
$$
f:=(\phi_1-\phi_0)+\delta w
$$
has a local maximum at $x_*$.
\item The function $\phi_1$ satisfies
the one--sided Hessian bound
\begin{equation}\label{eq:hess}
\scst I\le D^2\phi_1 (x_*) 
\end{equation}
\item If $U\subset\Re^n\times\Re\times\Re^n$ is the set in
Definition~\ref{uniform} then $\phi_0$ is $U$--admissible and satisfies
$\mean[\phi_0](x_*)\le\ecst$.  As we are assuming that $H_0=0$ and from the
definition of being uniformly elliptic $|b(x,\phi_0,D\phi_0)|\le \ecst$,
this implies that at the point $x_*$
\begin{equation}\label{phi0-bd}
\sum_{i,j}a^{ij}(x_*,\phi_0,D\phi_0)D_{ij}\phi_0(x_*)  \le 2\ecst
\end{equation}
\end{enumerate}
~\qed
\end{definition}

\subsection{Some Calculus Lemmata}
\label{sec:lemmata}

\begin{lemma}\label{test-fcn}\sl
Let $w:=\|x\|^{-\alpha}$ where $\alpha>0$ and $3r_0\le 1$ are as in the
standard setup.  Then $w$ is a decreasing
function of $\|x\|$ and for $r_0\le \|x\|\le 3r_0$ the inequalities
\begin{equation}\label{w-grad}
\|Dw\|(x)\le \alpha\|x\|^{-(\alpha+2)} \le \alpha r_0^{-(\alpha+2)}
\end{equation}
and
\begin{equation}\label{w-Hess}
-\alpha r_0^{-(\alpha+2)}I\le D^2w\le \alpha(\alpha+2)r_0^{-(\alpha+2)}I
\end{equation}
hold.
\end{lemma}

\begin{proof} By direct calculation
\begin{equation}\label{w-diff}
\begin{split}
D_iw &=-\alpha\|x\|^{-(\alpha+2)}x^i,\\
D_{ij}w &=\alpha(\alpha+2)
 \|x\|^{-(\alpha+4)}x^ix^j-\alpha\|x\|^{-(\alpha+2)}\delta_{ij}.
\end{split}
\end{equation}
From this  it follows $\|Dw\|(x)=\alpha\|x\|^{-(\alpha+2)}\|x\|\le
\alpha\|x\|^{-(\alpha+2)} \le \alpha r_0^{-(\alpha+2)}$ as
$r_0\le\|x\| \le1$.
The matrix $[x^ix^j]$ satisfies $0\le [x^ix^j]\le \|x\|^2I$ which, along
with~(\ref{w-diff}), implies the Hessian bounds~(\ref{w-Hess})
\end{proof}

The following lemma is not used in the proof of the analytic maximum
principle, but is important in the proof of the geometric version of
the maximum principle.

\begin{lemma}\label{grad-diff} Let $w$ be as in the Lemma~\ref{test-fcn}. 
If $f=\phi_1-\phi_0+\delta w$ has a local
maximum at $x_*\in B(x_1,r_1)$ with $r_1\le r_0$ and $\|x_*\|\ge r_0$ then
$$
\|D\phi_1(x_*)-D\phi_0(x_*)\|\le \delta\alpha r_0^{-(\alpha+2)}.
$$
\end{lemma}
\begin{proof} At a local maximum  $Df(x_*)=0$.  Using the estimate on
$\|Dw\|$ from Lemma~\ref{test-fcn}, the equation $Df(x_*)=0$ implies
$\|D\phi_1(x_*)-D\phi_0(x_*)\|=\delta\|Dw(x_*)\|\le \delta\alpha
r_0^{-(\alpha+2)}$.
\end{proof}

\begin{lemma}\label{lemma:matrix}
Let $A$ and $B$ be $n\times n$ symmetric matrices with $A$ positive definite
and let $c_1$, $c_2$, $c_3$ and $c_4$ be positive constants so that
$$
\frac1{c_1}I\le A\le c_2I,\quad B\ge -c_3  I, \quad \trace(AB)\le c_4
$$
then
$$
B\le c_1\bigl((n-1)c_2c_3+c_4\bigr)I.
$$
\end{lemma}
\begin{proof} By doing an orthogonal change of basis we can assume $B$ is
diagonal, say $B=\operatorname{Diag}(\beta^1\cd \beta^n)$.  Let $A=[a^{ij}]$.
Then
$$
\trace(AB)=a^{11}\beta^1+a^{22}\beta^2+\dots+a^{nn}\beta^n \le c_4.
$$
Solving this inequality for $\beta^1$ and using that the bounds on $A$
and $B$ imply the diagonal elements of $A$ satisfy $(1/c_1)\le a^{ii}\le
c_2$ and
$\beta^i\ge -c_3$, so that in particular
$-a^{ii}\beta^i\le c_2c_3$ and $1/a^{11}\le c_1$.  Thus
\begin{align*}
\beta^1&\le \frac{1}{a^{11}}(-a^{22}\beta^2-\dots-a^{nn}\beta^n+c_4)\\
&\le \frac{1}{a^{11}}\bigl((n-1)c_2c_3+c_4\bigr)\\
&\le c_1\bigl((n-1)c_2c_3+c_4\bigr),
\end{align*}
and a similar calculation shows $\beta^i\le c_1\((n-1)c_2c_3+c_4\)$ for all
$i$.  But as the $\beta^i$ are the eigenvalues of $B$ this implies
$B\le\((n-1)c_2c_3+c_4\)I$ and completes the proof.
\end{proof}

\begin{lemma}\label{Hess-bd}  Let $x_0$, $x_1$, $x_*$, $r_0$, $r_1$,
$\alpha$, $\delta$, $\phi_0$, $\phi_1$ be as in the standard
setup. Then
\begin{multline}\label{sum-bd}
|D^2\phi_0|(x_*)+|D^2\phi_1|(x_*) \le \\ 2\( \ecst^2\bigl((n-1)\bigl(\scst
+\delta\alpha r_0^{-(\alpha+2)}\bigr)
+2\bigr) +\delta\alpha r_0^{-(\alpha+2)} \). 
\end{multline}
\end{lemma}

\begin{proof}
The Hessian of $f$ is
$ D^2f=(D^2\phi_1-D^2\phi_0) + \delta D^2w $. Part of the standard setup
is that the  function  $f=(\phi_1-\phi_0)+\delta w$ has a local maximum
at $x_*$, and at a local maximum the Hessian satisfies $D^2f\le0$. Therefore
at  $x_*$  we have
$D^2\phi_1-D^2\phi_0  +\delta D^2w\le 0$.
Solving for $D^2\phi_0$ in this inequality and  using the inequality
$-D^2w(x)\le \alpha \|x\|^{-(\alpha+2)} I$ from
Lemma~\ref{test-fcn} and the inequality $D^2\phi_1(x_*)\ge -\scst I$ from
the standard setup we find that at $x_*$,
$$
D^2\phi_0(x_*) \ge D^2\phi_1(x_*)+\delta D^2w(x_*)
        \ge-\(\scst +\delta \alpha r_0^{-(\alpha+2)}\)I.
$$
This gives a lower bound on $D^2\phi_0(x_*)$.  Define symmetric matrices
$A:=[a^{ij}(x_*,\phi_0,D\phi_0)]$ and $B=D^2\phi_0(x_*)$.  Then as $\mean$
is uniformly elliptic the inequalities $(1/\ecst)I\le A\le \ecst I$ hold
and we have just derived the lower bound
$B\ge -\(\scst +\delta \alpha r_0^{-(\alpha+2)}\)I$. By  
inequality~(\ref{phi0-bd}) of the standard set up the inequality
$\trace(AB)\le 2\ecst$ holds.  Thus by the last lemma,
\begin{align*}
D^2\phi_0(x_*)&\le \ecst\((n-1)\ecst\bigl(\scst + \delta\alpha
r_0^{-(\alpha+2)}\bigr)+2\ecst\)I\\
&=\ecst^2\((n-1)\bigl(\scst+\delta\alpha r_0^{-(\alpha+2)}\bigr)+2\)I.
\end{align*}
Again we use that $D^2f(x_*)=D^2\phi_1(x_*)-D^2\phi_0(x_*)+\delta
D^2w(x_*)\le 0$, but this time solve for $D^2\phi_1(x_*)$.
Using the upper bound we have just derived for $D^2\phi_0(x_*)$ and the bound
$-D^2w\le \alpha r_0^{-(\alpha+2)}I$ of inequality~(\ref{w-Hess}) we
obtain,
$$
D^2\phi_1(x_*) \le D^2\phi_0(x_*)-\delta D^2w(x_*) 
\le
\beta I , 
$$
where 
\begin{equation}\label{eq:betadef}
\beta = 
\ecst^2\bigl((n-1)\bigl(\scst
+\delta\alpha r_0^{-(\alpha+2)}\bigr)+2\bigr)+\delta\alpha
r_0^{-(\alpha+2)} .
\end{equation}

We now have upper and lower bounds on both $D^2\phi_0(x_*)$ and
$D^2\phi_1(x_*)$.  The largest of the constants to appear in these bounds
is the constant on the right hand side of the upper bound on
$D^2\phi_1(x_*)$. Therefore,
$$
- \beta I \le D^2\phi_0(x_*), D^2\phi_1(x_*)  \le \beta I  ,
$$
with $\beta$ given by (\ref{eq:betadef}.
But if a symmetric matrix $S=[S^{ij}]$  satisfies $-\beta I\le
S \le \beta I$ the entries satisfy $|S^{ij}|\le \beta$.  The
inequality~(\ref{sum-bd}) now follows.
\end{proof}

\begin{cor}\sl\label{Hess-bd2}  
Assume in addition to the hypotheses of the last lemma that $\delta$
satisfies 
\begin{equation}\label{r-delta}
\delta\le \ol{\delta}(\alpha):=\frac{r_0^{\alpha+2}}{\alpha}.
\end{equation}
Then
$$
|D^2\phi_0|(x_*)+|D^2\phi_1|(x_*)\le \hcst,
$$
where
\begin{equation}\label{Hess-cst}
\hcst:=2\(\ecst^2\bigl((n-1)\bigl(\scst +1\bigr)+2\bigr)+1\)
\end{equation}
only depends on $\ecst$, $\scst$ and the dimension~$n$.
\end{cor}

\begin{proof} The bound~(\ref{r-delta}) implies
$\delta\alpha r_0^{-(\alpha+2)}\le1$
and hence the result follows from~(\ref{sum-bd}).
\end{proof}

\begin{lemma}\sl\label{lemma:coeff} Let $\mean$ and $U$ be as in the
statement of Theorem~\ref{thm:max}.
Let $\phi_0$ and $\phi_1$ be $U$--admissible
$C^2$ functions defined on some open subset of $\Omega$.  Then
\begin{align*}
\mean[\phi_1]-\mean[\phi_0]
=&\sum_{i,j}A^{ij}[\phi_0,\phi_1]D_{ij}(\phi_1-\phi_0)
  +\sum_iB^i[\phi_0,\phi_1]D_i(\phi_1-\phi_0)\\
 &\qquad +C[\phi_0,\phi_1](\phi_1-\phi_0)
\end{align*}
where the coefficients $A^{ij}[\phi_0,\phi_1]=A^{ji}[\phi_0,\phi_1]$ satisfy
the estimates
\begin{equation}\label{A-est}
\frac{1}{\ecst}\|\xi\|^2\le
\sum_{i,j}A^{ij}[\phi_0,\phi_1]\xi_i\xi_j\le\ecst\|\xi\|^2,
\end{equation}
so that
\begin{equation}\label{A-est2}
 \frac1{\ecst}\le
A^{ii}[\phi_0,\phi_1]\le \ecst, \quad\mbox{and}\quad |A^{ij}[\phi_0,\phi_1]|\le
\ecst.
\end{equation}
The coefficients $B^i[\phi_0,\phi_1]$, and $C[\phi_0,\phi_1]$ satisfy
\begin{equation}\label{BC-est}
|B^i[\phi_0,\phi_1]|,\  |C[\phi_0,\phi_1]|
        \le n^2\ecst(| D^2\phi_0| +| D^2\phi_1|+1).
\end{equation}
\end{lemma}

\begin{proof} One of the hypotheses of Theorem~\ref{thm:max} is
that the fibers $U_x$ are convex.  Thus if both $\phi_0$ and $\phi_1$
are $U$--admissible so is $\phi_t:=(1-t)\phi_0+t\phi_1$ for $0\le t\le1$.
Therefore we can compute
\begin{align}
\mean[\phi_1]-\mean[\phi_0]
&=\int_0^1\frac{d}{dt}\mean[\phi_t]\,dt\nn\\
&=\sum_{i,j}\int a^{ij}(x,\phi_t, D \phi_t)\,dt\, D_{ij}(\phi_1-\phi_0)\nn\\
&\qquad+\sum_{i,j,k}\int_0^1
\nn  \frac{\f a^{ij}}{\f p^k}(x,\phi_t, D \phi_t)
D_{ij}\phi_t\,dt\,D_k(\phi_1-\phi_0)\\
&\nn\qquad\qquad +\sum_i\int_0^1\frac{\f b}{\f p^i}(x,\phi_t, D
\phi_t)\,dt\,  D_i(\phi_1-\phi_0)\nn\\
&\nn\qquad+ \sum_{i,j}\int_0^1\frac{\f a^{ij}}{\f r}(x,\phi_t, D
\phi_t) D_{ij}\phi_t\,dt\, (\phi_1-\phi_0)\\
&\nn\qquad\qquad        +\int_0^1\frac{\f b}{\f r}(x,\phi_t, D
\phi_t)\,dt\,(\phi_1-\phi_0)\nn\\
\label{eq:diff}&=\sum_{i,j}A^{ij}[\phi_0,\phi_1] D_{ij}(\phi_1-\phi_0)
   +\sum_iB^i[\phi_0,\phi_1] D_i(\phi_1-\phi_0)\\
&\qquad\nn+C[\phi_0,\phi_1](\phi_1-\phi_0),
\end{align}
where, by re-indexing some of the sums involved, we find the coefficients
$A^{ij}$, $B^i$ and $C$ have the formulas
\begin{align}\label{A-def}
A^{ij}[\phi_0,\phi_1]&:=\int_0^1a^{ij}(x,\phi_t, D \phi_t)\, dt,\\
\label{B-def}
B^i[\phi_0,\phi_1]&:=
        \int_0^1\biggl(\sum_{j,k}\frac{\f a^{jk}}{\f p^i}(x,\phi_t, D
\phi_t) D_{jk}\phi_t+
        \frac{\f b}{\f p^i}(x,\phi_t, D \phi_t)\biggr) dt,\\
\label{C-def}
C[\phi_0,\phi_1]&:=\int_0^1\biggl(\sum_{i,j}\frac{\f a^{ij}}{\f r}(x,\phi_t, D
\phi_t) D_{ij}\phi_t+\frac{\f b}{\f r}(x,\phi_t, D \phi_t)\biggr) dt.
\end{align}
Then~(\ref{A-est}) follows
from~(\ref{a-est}) by integration.  The inequalities~(\ref{A-est2}) are
algebraic consequences of~(\ref{A-est}).
From $\phi=(1-t)\phi_0+t\phi_1$ for $t\in[0,1]$ we get,
$$
| D^2\phi_t|\le
(1-t)| D^2\phi_0|+t| D^2\phi_1|
\le | D^2\phi_0|+| D^2\phi_1|,
$$
so that if $\ecst$ is the constant in~(\ref{der-bd}),
$$
\begin{aligned}
\left|\int_0^1\frac{\f a^{jk}}{\f p^i}(x,\phi_t, D \phi_t)
D_{jk}\phi_t\,dt\right| &\le \ecst(| D^2\phi_0|+| D^2\phi_1|), \\
\left|\int_0^1\frac{\f b }{\f p^i}(x,\phi_t, D \phi_t) \,dt\right|
&\le
\ecst.
\end{aligned}
$$
Using this in the formula~(\ref{B-def}) defining $B^i[\phi_0,\phi_1]$ and the
inequality $1\le n^2$ we obtain,
$$
|B^i[\phi_0,\phi_1]|\le
\ecst(n^2(| D^2\phi_0|+| D^2\phi_1|)+1)
\le n^2\ecst(| D^2\phi_0|+| D^2\phi_1|+1)
$$
as required.  The derivation of the bound on $C[\phi_0,\phi_1]$ is
identical.
\end{proof}

\begin{cor}\sl\label{coeff-bd}  Let $x_0$, $x_1$, $x_*$, $r_0$, $r_1$,
$\alpha$, $\delta$, $\phi_0$, $\phi_1$ be as in the
standard setup and assume  $\delta\le
\ol{\delta}(\alpha)$, as given in~(\ref{r-delta}).  Then
$$
|B^i[\phi_0,\phi_1]|(x_*),\ |C[\phi_0,\phi_1]|(x_*)\le n^2\ecst(\hcst+1),
$$
and
\begin{equation}\label{norm-bar}
\|B[\phi_0,\phi_1]\|\le n^3\ecst(\hcst+1),
\end{equation}
where $\hcst$ is as in equation~(\ref{Hess-cst}).
\end{cor}

\begin{proof} This follows by using the estimates~(\ref{BC-est})
and~(\ref{sum-bd}) in the last lemma, and in the estimate for
$\|B[\phi_0,\phi_1]\|$ the inequality $n^{5/2}\le n^3$ was used.
\end{proof}

\begin{lemma}\sl \label{op-bd}
Let $x_0$, $x_1$, $x_*$, $r_0$, $r_1$, $\alpha$, $\delta$,
$\phi_0$, $\phi_1$ be as in the standard setup and set
\begin{equation}\label{eq:alpha-bar}
\alpha=\ol{\alpha}:=-2+\ecst\(1+n\ecst+n^3\ecst(\hcst+1)\)
\end{equation}
(This makes $\alpha$ a solution to
$\((\alpha+2)/\ecst-n\ecst- n^3\ecst(\hcst+1)\)=1$). Also
assume $\delta$ satisfies the
inequality~(\ref{r-delta}), and let $L_{\phi_0,\phi_1}$ be
as defined in Equation~(\ref{Lbar}).  Then
$$
L_{\phi_0,\phi_1}w(x_*)\ge 1.
$$
\end{lemma}

\begin{proof} We first note for future reference that $n,\ecst\ge1$, $\hcst>0$
implies $\ol{\alpha}\ge -2+3=1$.
Under the hypotheses of the lemma, the bounds in~(\ref{A-est2}) on
$A^{ij}[\phi_0,\phi_1]$ and
the bound~(\ref{norm-bar}) on $\left\|B[\phi_0,\phi_1]\right\|(x_*)$ hold.
We will also use the bound $\|Dw\|(x)\le\alpha\|x\|^{-(\alpha+2)}$ of
Lemma~\ref{test-fcn}.  To simplify notation
we assume  all functions are evaluated at $x=x_*$ and
drop the subscript of~$*$.
Using (\ref{w-grad}) and the explicit formula (\ref{w-diff}) for
$D_{ij}w$ we have,
\begin{align*}
L_{\phi_0,\phi_1}w(x_*)&=
\sum_{i,j}A^{ij}[\phi_0,\phi_1]D_{ij}w+\sum_i{B}^i[\phi_0,\phi_1]D_iw\\
&=\alpha(\alpha+2)\|x\|^{-(\alpha+4)}\sum_{i,j}A^{ij}[\phi_0,\phi_1]x^ix^j
         \\
 & \hskip 0.3truein -\alpha\|x\|^{-(\alpha+2)}\sum_iA^{ii}[\phi_0,\phi_1] 
+\sum_i{B}^i[\phi_0,\phi_1]D_iw\\
&\ge\alpha(\alpha+2)\|x\|^{-(\alpha+4)}\dfrac{\|x\|^2}{\ecst}
        -\alpha\|x\|^{-(\alpha+2)}
        n\ecst-\left\|{B}[\phi_0,\phi_1]\right\|\|Dw\|\\
&\ge\alpha\|x\|^{-(\alpha+2)}\(\dfrac{\alpha+2}{\ecst}-n\ecst\)
        -n^3\ecst(\hcst+1)\alpha\|x\|^{-(\alpha+2)}\\
&=\alpha\|x\|^{-(\alpha+2)}\(\dfrac{\alpha+2}{\ecst}
                -n\ecst-n^3\ecst(\hcst+1)\)\\
&=\alpha\|x\|^{-(\alpha+2)}\ \ge\  1,
\end{align*}
where the last couple of steps hold because of the choice of
$\alpha=\ol{\alpha}$ (so that also $\alpha\ge1$), and
$\|x\|^{-(\alpha+2)}\ge1$ as $\|x\|\le 3r_0\le 1$.  This completes
the proof.
\end{proof}

\subsection{Proof of the Maximum Principle}\label{proof}

We follow the outline given in \S \ref{sec:outline} and assume
that we have already made the normalizations given in \S \ref{setup}.
We then choose $\alpha$ as in Lemma~\ref{op-bd}:
$$
\alpha=\ol{\alpha}=- 2+\ecst\(1+n\ecst+n^3\ecst(\hcst+1)\)
$$
This is a formula for $\alpha$ in terms of just the
parameters $\ecst$, $\scst$ and the dimension~$n$ (recall the definition of
$\hcst$ given in equation~(\ref{Hess-cst})).  Let
$$
r_1=\min\left\{r_0,\
        \frac1{4\(n^2\ecst(\hcst+1)\)\(\alpha r_0^{-(\alpha+2)}\)}\right\}
$$
Thus, keeping in mind the definition of $\alpha$,  $r_1$ only depends on $r_0$,
$\ecst$, $\scst$, and $n$. Let $S'$, $S''$ as in \S \ref{sec:outline} and
let $h:=(u_1-u_0) + \delta (w-w(x_1))$ where we now
determine how to choose $\delta$.
The function $u_1-u_0$ is negative and upper semi-continuous
on the compact set $S'$ and so there is a $\delta_1>0$ such that
$ (u_1-u_0) + \delta_1 (w -w(x_1)) <0$ on $S'$.  Let
$$
\delta:=\min\{\delta_1, \ol{\delta}(\alpha)\},
$$
where $\ol{\delta}(\alpha)$ is given by~(\ref{r-delta}) (so $\delta$
is defined just in terms of $\ecst$, $\scst$, $n$, $r_0$, and
$\delta_1$).  As $\delta\le\delta_1$ we have $h<0$ on $S'$.  The
function $w$ is a decreasing function of $\|x-x_0\|=\|x\|$ so the
function $h$ will be negative on the set $S''$.  Therefore $h$ is
negative on all of $\f B(x_1,r_1)=S'\cup S''$.  But $x_1$ was chosen
so that $u_1(x_1)=u_0(x_1)$ and thus
$h(x_1)=(u_1(x_1)-u_0(x_1))+\delta(w(x_1)-w(x_1))=0$. Hence  $h$ has
an local maximum at some interior point $x_*$ of $B(x_1,r_1)$.

Choose $\e>0$ so that
$$
\e <\min\left\{\frac{\delta}{4}, \ecst\right\}.
$$
Recall that  we are assuming  $H_0=0$. Let $\phi_0$
be an upper support function for $u_0$ at
$x_*$ and $\phi_1$ a lower support function for $u_1$ at $x_*$ so that
$\mean[\phi_0]\le \e$, $\mean[\phi_1]\ge -\e$ and the one--sided Hessian
bound~(\ref{eq:hess}) holds.   Using that
$\phi_0$ and $\phi_1$ are upper and lower support functions at $x_*$
we see  $f:= (\phi_1-\phi_0)+\delta w$  has a
local maximum at~$x_*$.  Therefore $x_0$, $x_1$, $x_*$, $r_0$, $r_1$,
$\alpha$, $\delta$, $\phi_0$, $\phi_1$ satisfy all the
conditions of the standard setup together with the
hypotheses of Lemma~\ref{op-bd}.  Using
formula~(\ref{com}), the definition~(\ref{Lbar}) of
$L_{\phi_0,\phi_1}$, and the choice of $\e$
\begin{align*}
-\frac{\delta}{2}&\le -2\e  \le \mean[\phi_1](x_*)-\mean[\phi_0](x_*)\\
&= L_{\phi_0,\phi_1}\(\phi_1-\phi_0\)(x_*)
        +{C}[\phi_0,\phi_1]\(\phi_1-\phi_0\)(x_*),
\end{align*}
so that
$$
 L_{\phi_0,\phi_1}\(\phi_1-\phi_0\)(x_*)\ge -\frac{\delta}{2}
                        -|C[\phi_0,\phi_1]||\phi_1(x_*)-\phi_0(x_*)|.
$$
By Lemma~\ref{op-bd} $L_{\phi_0,\phi_1}w(x_*)\ge 1$.
Thus for $f=(\phi_1-\phi_0)+\delta w$ we use the last displayed
inequality, the bound of Corollary~\ref{coeff-bd} on $|C[\phi_0,\phi_1]|$,
and the equality  $u_i(x_*)=\phi_i(x_*)$ to  compute
\begin{align}
L_{\phi_0,\phi_1}f(x_*)&\ge L_{\phi_0,\phi_1}(\phi_1-\phi_0)(x_*)+\delta
 \ge \frac\delta2-|C[\phi_0,\phi_1]|(\phi_1(x_*)-\phi_0(x_*)|\nn\\
\label{pre-contra} &\ge \frac\delta2-n^2\ecst(\hcst+1)|u_1(x_*)-u_0(x_*)|.
\end{align}
We now derive an estimate on $|u_1(x_*)-u_0(x_*)|$.  From the hypotheses of
the theorem $u_1(x_*)-u_0(x_*)\le 0$.  The function $h=(u_1-u_0)+\delta
(w-w(x_1))$ has its maximum in the ball $B(x_1,r_1)$ at the point $x_*$ and
$h(x_1)=0$ and thus $h(x_*)=(u_1(x_*)-u_0(x_*))+\delta(w(x_*)-w(x_1))\ge0$.
Therefore, using the bound~(\ref{w-grad}) on $\|Dw\|$, we have
$$
|u_1(x_*)-u_0(x_*)|\le \delta |w(x_*)-w(x_1)|\le \delta \|x_*-x_1\|\alpha
r_0^{-(\alpha+2)}\le \delta \alpha r_0^{-(\alpha+2)} r_1 .
$$
Using this in~(\ref{pre-contra}) along with the definition of $r_1$ gives
$$
L_{\phi_0,\phi_1}f(x_*)\ge \frac\delta2-\frac\delta4=\frac\delta4>0.
$$
But from the first and second derivative tests of calculus $Df(x_*)=0$
and $D^2f(x_*)\le 0$ so
$$
L_{\phi_0,\phi_1}f(x_*)
        =\sum_{i,j}A^{ij}[\phi_0,\phi_1]D_{ij}f(x_*)\le0.
$$
This contradicts $L_{\phi_0,\phi_1}f(x_*)>0$ and completes the proof
$u_0\equiv u_1$.
\bigskip

We now set $u:=u_0=u_1$ and show $u$ is locally a $C^{1,1}$ function,
i.e. that  $Du(x)$ exists for all $x_0\in\Omega$ and it locally satisfies
a Lipschitz condition $\|Du(x_1)-Du(x_0)\|\le C\|x_1-x_0\|$ for some
$C\ge0$.  First, as $u=u_1=u_0$ is both upper and lower semi-continuous
it 
is continuous.  
As $u_0=u_1$ the function $u_1-u_0\equiv0$ has a
local maximum at every point of its domain $\Omega$. 
Letting $x$ be any point in $\Omega$, 
let $\phi_0$ be a $C^2$ upper support
function for $u_0$ at $x$ with $\mean[\phi_0]\le\ecst$ at $x$ and
let $\phi_1$ be a $C^2$ lower support
for $u_1$ at $x$ which satisfies
the one--sided Hessian bound $D^2\phi_1(x) \ge - \scst I$.  
Then $\phi_1 - \phi_0$ has a local maximum at $x$ and a simplified
version of the proof of Lemma~\ref{Hess-bd} (with $\delta =0$) shows there 
is a constant $\hcst'$ independent of $x$ so that 
$-\hcst' I\le D^2\phi_0, D^2\phi_1\le \hcst' I$.  (Going through the proof of
Lemma~\ref{Hess-bd} with $\delta=0$ shows that
$\hcst'=\ecst^2\bigl((n-1)\scst+2)$ works.)

The following is certainly known, but as we have not found an explicit
reference we include a short proof.

\begin{lemma}\sl\label{global-spt}
Let $v$ be a continuous function on the convex open set $\Omega$ and assume
that $v$ satisfies $D^2v\ge -CI$
in the sense of support functions for
some $C\in \Re$. (That is for every $x\in \Omega$ there is a $C^2$ lower support
function $\phi$ to $v$ at $x$ so that $D^2\phi\ge -  CI$ near $x$.)  Then
for each $x_0\in \Omega$ there is a vector $a\in\Re^n$ so that for all $x\in
\Omega$
$$
v(x)\ge v(x_0)+\la x-x_0,a\ra -\frac12 C\|x-x_0\|^2.
$$
If $v$ is differentiable at $x_0$ then $a=Dv(x_0)$.
\end{lemma}

\begin{proof} We first give the proof in the one dimensional case.  Then
$\Omega\subseteq \Re$ is an interval.  Let $\phi$ be a lower support
function to $v$ at $x_0$ that satisfies $D^2\phi=\phi''\ge -C$ near
$x_0$ and let $a=D\phi(x_0)$.  We claim $v(x)\ge
v(x_0)+a(x-x_0)-(C/2)(x-x_0)^2$.  Let
$f(x):=v(x_0)+a(x-x_0)-(C/2)(x-x_0)^2$.  As $\phi(x_0)=f(x_0)$,
$\phi'(x_0)=f'(x_0)$ and $\phi''(x)\ge -C=f''(x)$ we have an interval
about $x_0$ so that $v(x)\ge \phi(x)\ge f(x)$.  Assume, toward a
contradiction, there is an $x_1$ so that $v(x_1)-f(x_1)<0$.  Then the
function $v_1:=v-f$ will satisfy $v_1''\ge 0$ in the sense of support
functions and $v_1(x_0)=0$, $v_1(x_1)\le 0$ and the function
$v_1(x)\ge 0$ for $x$ near $x_0$ so the function $v_1$ will have an
maximum at a point $x_*$ between $x_0$ and $x_0$.  Let $v_0$ be the
constant function $v_0(x)=v_1(x_*)$.  Then $v_1\le v_0$,
$v_1(x_*)-v_0(x_*)$, and $v_1''\ge 0$ (in the sense of support
functions) and $v_0''=0$ (in the strong sense) so the one variable
case of of the linear maximum principle (which is easy to verify)
implies $v_0(x)=v_1(x)$ for $x$ between $x_0$ and $x_1$.  As
$v_1(x_0)=0$ and $v_0$ is constant this implies $v_1(x)=v(x)-f(x)$ for
all $x$ between $x_0$ and $x_1$.  This in particular implies
$0=v(x_1)=f(x_1)$ which contradicts the assumption $v(x_1)-f(x_1)<0$
and completes the proof in the one dimensional case.

We  return to the general case.  Let $\phi$ be a lower support
function for $v$ at $x_0$ that satisfies $D^2\phi\ge -CI$ near $x_0$
and let $a:=D\phi(x_0)$.  For any unit vector $b\in\Re^n$ let
$v_b(t)=v(x_0+tb)- t \la b,a\ra$.  Then a lower support $\psi$
function to $v$ at $x_0+t_0b$ that satisfies $D^2\psi\ge -CI$ yields
the lower support function $\psi_b(t):=\psi(x_0+tb)$ to $v_b$ at $t_0$
that satisfies $\psi_b''(t)\ge -C$ near $t_0$.  Thus $v_b$ satisfies
the one variable version  of the result and so
$v_b(t)=v(x_0+tb)\ge v(x_0)-t \la b,a\ra -(C/2)t^2$.  But as $\Omega$ is
convex every point of $\Omega$ can be written as $x=x_0+tb$ for some $t$
and some unit vector and so the multidimensional case reduces to the one
dimensional case.  This completes the proof.
\end{proof}

The last lemma implies at each point of $\Omega$ that $u$ has global
upper and lower support paraboloids with ``opening'' $2\hcst$ (that
is, in the terminology of Caffarelli and
Cabre~\cite{Caffarelli-Cabre}, the second order term of the paraboloid
is $\pm\hcst\|x\|^2$).  It now follows that $u$ is of class $C^{1,1}$, see
for example~\cite[Prop~1.1~p7]{Caffarelli-Cabre}.

Now assume $a^{ij}$ and $b$ are locally of class $C^{0,\alpha}$ of all its
arguments $(x,r,p)$ for some $\alpha\in
(0,1)$ (this is in addition to the assumption that they are $C^1$ functions
of the arguments $(r,p)$) and we will show $u=u_0=u_1$ is locally
$C^{2,\alpha}$.  Define a linear second order  differential operator by
$$
{\mathcal A}_uf:=\sum_{i,j}a^{ij}(x,u(x),Du(x))D_{ij}f.
$$
As $Du$ is locally Lipschitz the functions $x\mapsto
a^{ij}(x,u(x),Du(x))$ and $x\mapsto b(x,u(x),Du(x))$  are locally of
class $C^{0,\alpha}$.  Let $B$ be an open ball whose closure is contained
in $\Omega$.  Then by a standard existence result
of the linear Schauder theory (cf.~\cite[Thm~6.13~p101]{Gilbarg-Trudinger})
gives there is a unique function $v$ continuous on the closed ball
$\ol{B}$, locally of class $C^{2,\alpha}$ in $B$ and so that $v$ solves the
boundary value problem
$$
{\mathcal A}_uv(x)=-b(x,u(x),Du(x))+H_0\quad\mbox{in $B$},\qquad v=u\quad\mbox{on
        $\f B$}.
$$
But $u$ solves the same boundary value problem
$$
{\mathcal A}_uu(x)=-b(x,u(x),Du(x))+H_0\quad\mbox{in $B$},\qquad u=u\quad\mbox{on
        $\f B$}
$$
but in the sense of support functions rather than in the classical sense.
Thus the function $f:=u-v$ will satisfy ${\mathcal A}_uf=0$ in $B$ in the
sense of support functions and $f=0$ on $\f B$.  Then $f=0$ in $B$ by
Calabi's version of the Hopf maximum principle~\cite{Calabi:maximum}.
Therefore $u=v$ in $B$ so that $u$ is locally of class~$C^{2,\alpha}$ in
$B$.  As every point of $\Omega$ is contained in such a ball $B$ the
function $u$ is locally of class $C^{2,\alpha}$.

We use this as the base of an induction to prove higher regularity.
Assume for some $k\ge1$ the functions $a^{ij}$ and $b$ are locally
$C^{k,\alpha}$ and $u$ is locally  $C^{k+1,\alpha}$.  Then
$x\mapsto a^{ij}(x,u(x),Du(x))$ and $x\mapsto b(x,u(x),Du(x))$ are locally
$C^{k,\alpha}$. But $u$ is a solution to the linear
equation ${\mathcal A}_uu(x)=-b(x,u(x),Du(x))+H_0$ and by the linear Schauder
regularity theory (cf.~\cite[Thm~6.17~p104]{Gilbarg-Trudinger}) this implies
$u$ is locally of class $C^{k+2,\alpha}$.  This completes the induction
step and finishes the proof of Theorem~\ref{thm:max}.

\section{Geometric Maximum Principles for Hypersurfaces in  Lorentzian
and Riemannian Manifolds}
\label{sec:geomax}


The version of the analytic maximum principle given by
Theorem~\ref{thm:max} is especially natural in the Lorentzian setting
as $C^0$~spacelike hypersurfaces (Definition~\ref{def:spacelike}
below) can always be locally represented as graphs.
Theorem~\ref{thm:max} also applies to hypersurfaces in Riemannian
manifolds that can be represented locally as graphs (see
Theorem~\ref{Riem-max}).  However in the Riemannian setting there are
many cases, such as horospheres, where one wants to use the maximum
principle but this graph condition does not {\em a priori\/} hold.  A
version of the geometric maximum principle adapted to horospheres and
other rough hypersurfaces is given in~\cite{Howard:inner-sphere}.  We
now give a detailed proof of the geometric maximum principle for
$C^0$spacelike hypersurfaces in Lorentzian manifolds.

%

We first state our conventions on the sign of the second fundamental
form and the mean curvature.  To fix our choice of
signs, a Lorentzian manifold $(M,g)$ is an $n$-dimensional manifold
$M$ that has semi-Riemannian~$g$ metric with signature $+,\dots,+,-$.
Our results are basically local and, since every Lorentzian manifold is
locally time orientable, we assume that all our Lorentzian manifolds
are timed oriented and use the standard terminology of {\bi
spacetime}\/ for a time oriented Lorentzian manifold.
See \cite{Beem-Ehrlich:Lorentz2} for background
on Lorentzian geometry.

Given a smooth spacelike hypersurface $N\subset M$ we will always use the
future pointing unit normal $\nor$ which, as $N$ is spacelike, will be
timelike, $g(\nor,\nor)=-1$.  
Let $\nabla$ be the metric connection on $(M,g)$.  Then the
second fundamental form $h$ and mean curvature $H$ of $N$ are defined
by
$$
h(X,Y):=-g(\nabla_XY,\nor),\qquad H:=\frac1{n-1}\trace_{g\big{|}_{N}}\!\! h
=\frac1{n-1}\sum_{i=1}^{n-1}h(e_i,e_i)
$$
where $X$, $Y$ are smooth vectors fields tangent to $N$ and $e_1,\cd
e_{n-1}$ is a smooth locally defined orthonormal frame field along
$N$.  

If $(M,g)$ is flat Minkowski space with metric 
$g:=(dx^1)^2+\cdots+(dx^{n-1})^2-(dx^n)^2$ and if $N$ is given by a
graph $x^n=f(x^1\cd x^{n-1})$ with $\|D f\|^2:=(\f f/\f
x^1)^2+\cdots (\f f/\f x^{n-1})^2<1$ then with this choice of signs, 
the mean curvature of $N$ is given by
$$
H=\frac1{n-1}\sum_{i=1}^{n-1}\frac{\f}{\f x^i}\(\frac1{\sqrt{1-\|D f
\|^2}}\frac{\f f}{\f x^i}\).
$$
(Here $\f/\f x^n$ is in the direction of the future.)
Thus the symbol of the linearization of $H$ is positive definite.
This is also the choice of sign so that if $N$ is moved along the future
pointing normal $\nor$ then $N$ is expanding when $H$ is positive and
shrinking when $H$ is negative, where expanding and shrinking are
measured in terms of expansion or shrinking of the area element of $N$.
Thus very roughly $H$ is the local ``Hubble constant'' of an observer
on $N$ who believes that the world-lines of rest particles are the
geodesics normal to $N$.
\bigskip

In the maximum principle we wish to use the weakest natural notion of a
hypersurface being spacelike.  The following definition was introduced
in~\cite{Eschenburg-Galloway} or
see~\cite[p.~539]{Beem-Ehrlich:Lorentz2}. 

\begin{definition}\label{def:spacelike}
A subset $N\subset M$ of the spacetime $(M,g)$ is a {\bi $C^0$ spacelike
hypersurface}\/ iff for each $p\in N$ there is a neighborhood $U$ of $p$ in
$M$ so that $N\cap U$ is acausal and edgeless in $U$.
~\qed
\end{definition}

\begin{remark}\label{local-gl-hy}
In this definition note that if $D(N\cap U,U)$ is the domain of
dependence of $N\cap U$ in $U$ then $D(N\cap U,U)$ is open in $M$ and
$N\cap U$ is a Cauchy hypersurface in $D(N\cap U,U)$.  But a spacetime
that has a Cauchy hypersurface is globally hyperbolic.  Thus by
replacing $U$ by $D(N\cap U,U)$ we can assume the neighborhood $U$ in
the last definition is globally hyperbolic and that $N\cap U$ is a
Cauchy surface in $U$.  In particular, a $C^0$ spacelike hypersurface is
a topological (in fact Lipschitz) submanifold of codimension one.~\qed
\end{remark}

Let $(M,g)$ be a spacetime and let $N_0$ and $N_1$ be two $C^0$ spacelike
hypersurfaces in $(M,g)$ which meet at a point $q$.  We say that
$N_0$ is {\bf locally to the future  of $N_1$ near $q$} iff for 
some neighborhood $U$
of $p$ in which $N_1$ is acausal and edgeless, $N_0 \cap U\subset
J^+(N_1,U)$,
where $J^+(N_1,U)$ is the causal future of $N_1$ in $U$.
In a time-dual fashion we may define what it means for $N_0$ to be locally
to the
past of $N_1$ near $q$.

Now consider a $C^0$  spacelike hypersurface $N$ in
a spacetime $(M,g)$. In the context of the following definition, $S$ is a
{\bi future support hypersurface}\/ for $N$ at $x_0 \in N$ iff $x_0\in S$
and $S$
is locally to the future of $N$ near $x_0$.  Time-dually, $S$ is a
{\bi past support hypersurface}\/ for $N$ at $x_0 \in N$ iff $x_0\in S$ and $S$
is locally to the past of $N$ near $x_0$.



\begin{definition}\label{def:hyper-mean}
Let $N$ be a $C^0$ spacelike hypersurface in the spacetime $(M,g)$ and
$H_0$ a constant.  Then
\begin{enumerate}
\item $N$ has {\bi mean curvature $\le H_0$ in the sense of support
hypersurfaces}\/ iff for all $q\in N$ and $\e>0$ there is a $C^2$ future
support hypersurface $S_{q,\e}$ to $N$ at $q$ and the mean
curvature of $S_{q,\e}$ at $q$ satisfies
$$
H_q^{S_{q,\e}}\le H_0+\e.
$$
\item $N$ has {\bi mean curvature $\ge H_0$ in the sense of support
hypersurfaces with one--sided Hessian bounds}\/ iff for all compact
sets $K\subseteq N$ there is a compact set $\widehat{K}\subset T(M)$
and a constant $C_K >0$ such that for all $q\in K$ and $\e>0$ there is a
$C^2$ past support hypersurface $P_{q,\e}$ to $N$ so that
\begin{enumerate}
\item The future pointing unit normal $\nor^{P_{q,\e}}(q)$ to
$P_{q,\e}$ at $q$ is in $\widehat{K}$.
\item At the point $q$ the mean curvature $H^{P_{q,\e}}$ and second
fundamental form $h^{P_{q,\e}}$ of $P_{q,\e}$ satisfy
$$
H_q^{P_{q,\e}}\ge H_0-\e,\quad h^{P_{q,\e}}_q
\ge - C_K g\big{|}_{P_{q,\e}}.
$$
\end{enumerate}
\end{enumerate}
~\qed
\end{definition}

\begin{remark} \label{geo-prob}
As will be seen below (Lemma~\ref{cpt:vectors} and the discussion
following it) the condition that the unit normals to the
support hypersurfaces $P_{q,\e}$ to $N$ for $q\in K$ all remain in
a compact set $\widehat{K}$ is equivalent to the mean curvature
operator being uniformly elliptic on the set $K\subseteq N$.
Therefore the definition of $N$ having mean curvature $\ge H_0$ in the
sense of support hypersurfaces with one--sided Hessian bounds has
built into it that the mean curvature operator on $N$ is uniformly
elliptic.

However in the definition of $N$ having mean curvature $\le H_0$ in
the sense of support hypersurfaces there is no restriction on the
normals to the support hypersurfaces and so the mean curvature
operator need not be uniformly elliptic on $N$.  (In the set up
for the analytic maximum principle case this is equivalent to dropping
the assumption that the upper support functions $\phi_{x,\e}$ in
Definition~\ref{mean-spt} have to be $U$--admissible.)  Therefore in
proving the geometric version of the maximum principle we will have to
prove a (fortunately trivial) estimate that shows in fact the normals
to the future support hypersurfaces are well behaved.~\qed
\end{remark}

In practice it is often the case that the support hypersurfaces
$S_{q,\e}$ and $P_{q,\e}$ can be chosen to be geodesic spheres
through $q$ all tangent to each other and with increasing radii in
which case the one--sided bound on the second fundamental forms
of the past support hypersurfaces $P_{q,\e}$ in the definition of
mean curvature $\ge H_0$ in the sense of support functions with a one--sided
Hessian bound holds for easy {\em a priori}\/ reasons.  The
following proposition, whose proof just relies on the continuous
dependence of solutions of ordinary differential equations on initial
conditions and parameters, summarizes this.  The details are left to
the reader.

\begin{prop}\label{sphere-bds}
Let $(M,g)$ be a spacetime, $r_0>0$ and $K\subset T(M)$ a compact set of
future pointing timelike unit vectors.  Assume that there is a $\delta$ so
that for all $\eta \in K$, the geodesic 
$\gamma_\eta(t):=\exp(t\eta)$ maximizes the
Lorentzian distance on the interval $[0,r_0+\delta]$.  For each $\eta\in K$
and $r>0$ let $\pi(\eta)$ be the base point of $\eta$ and set
\begin{equation}\label{sphere-def}
S_{\eta,r}:=\{p: d(p,\exp(r\eta))=r\}.
\end{equation}
Then $S_{\eta,r_0}$ contains $\pi(\eta)$ and in a neighborhood of
$\pi(\eta)$ it is a smooth spacelike hypersurface whose future
pointing unit normal at $\pi(\eta)$ is $\eta$.  There is a uniform 
two--sided bound on the second fundamental forms
$h^{S_{\eta,r_0}}_{\pi(\eta)}$ (or equivalently, uniform
bounds on the absolute values of the principal curvatures) as $\eta$
ranges over $K$ (these bounds only depend on $(M,g)$, $K$, and $r_0$).
Thus if $\eta\in K$ and $S$ is a smooth spacelike hypersurface that
passes through $\pi(\eta)$ and locally near $\pi(\eta)$ is in the
causal future of $S_{\eta,r_0}$ then there is a lower bound on the
second fundamental form of $S$ at $\pi(\eta)$ that only depends on
$(M,g)$, $K$, and $r_0$.  In particular if $r\ge r_0$, and
$S_{\eta,r}$ is smooth near $\pi(\eta)$ then $S_{\eta,r}$ is in the
causal future of $S_{\eta,r_0}$ (by the reverse triangle inequality)
and so there is a lower bound on the second fundamental form of
$S_{\eta,r}$ only depending on $(M,g)$, $K$, and $r_0$.~\qed
\end{prop}

\begin{thm}[Lorentzian Geometric Maximum Principle]\label{geo:max}
Let $N_0$ and $N_1$ be $C^0$ spacelike hypersurfaces in a spacetime $(M,g)$
which meet at a point $q_0$, such that $N_0$ is locally to the future of
$N_1$ near
$q_0$.  Assume for some constant $H_0$:
\begin{enumerate}\item $N_0$ has mean curvature  $\le H_0$ in the
        sense of support hypersurfaces.
        \item $N_1$ has mean curvature $\ge H_0$ in the sense of
        support hypersurfaces with one--sided Hessian bounds.
\end{enumerate}
Then $N_0=N_1$ near $q_0$, i.e., there is a neighborhood $\mathcal O$
of $q_0$ such that
$N_0\cap {\mathcal O} = N_1 \cap {\mathcal O}$.  Moreover,
$N_0\cap {\mathcal O} = N_1 \cap {\mathcal O}$ is a
smooth  spacelike hypersurface with mean curvature $H_0$.
\end{thm}
In \S\ref{reduce} below, Theorem \ref{geo:max} is reduced to the analytic 
maximum principle. The proof is given in \S\ref{sec:proofgeomax}.



\begin{remark}\label{rem:finiteder}	
If the metric only has finite differentiability, say $g$ is
$C^{k,\alpha}$ with $k\ge 2$ and $0<\alpha<1$ then, since the functions
$a^{ij}$ and $b$ in the definition of the mean curvature operator
$\mathcal H$ (see the proof below) depend on the first derivatives of the
metric, they are of class $C^{k-1,\alpha}$. Thus the regularity
part of Theorem \ref{thm:max} implies the hypersurface 
$N_0\cap {\mathcal O} =N_1\cap {\mathcal O}$ 
in the statement of the last theorem is 
$C^{k+1,\alpha}$.~\qed
\end{remark}

Note that hypothesis 2 is trivially satisfied if $N_1$ is
a smooth spacelike
hypersurface with mean curvature $\ge H_0$ in the usual sense.
This yields a ``rough-smooth" version of our geometric maximum
principle which does not require
any a priori Hessian estimates.  This version implies the ad hoc
maximum principle for the level sets of the Lorentzian Busemann
function obtained in \cite{galloway:split} to prove the
Lorentzian splitting theorem (see \cite{galloway:warsaw} for further
applications).
However, the geometric maximum principle, when used in its full
generality, yields a more intrinsic
and conceptually simplified proof of the Lorentzian splitting theorem.
See the next section for further discussion of this point and additional
applications.


\subsection{Reduction to the Analytic Maximum Principle}\label{reduce}
Let $(M,g)$ be an $n$~dimensional spacetime and let $\nabla$ be the
metric connection of the metric $g$.  Then near any point $q$ of $M$
there is a coordinate system $(x^1,\dots,x^n)$ so that the metric
takes the form
\begin{equation}\label{eq:metric}
g=\sum_{A,B=1}^ng_{AB}dx^Adx^B=\sum_{i,j=1}^{n-1}g_{ij}dx^idx^j-(dx^n)^2
\end{equation}
and so that $\f/\f x^n$ is a future pointing timelike unit vector.
(To construct such coordinates choose any smooth spacelike
hypersurface $S$ in $M$ passing through $q$ and let $(x^1\cd x^{n-1})$
be local coordinate on $S$ centered on at $q$.  Let $x^n$ be the
signed Lorentzian distance from $S$.  Then near $q$ the coordinate
system $(x^1\cd x^n)$ is as required.)  
Let $f$ be a function defined near the origin in $\Re^{n-1}$ with
$f(0)=0$.  Then define a map $F_f$ from a neighborhood of the origin
in $\Re^{n-1}$ to $M$ so that in the coordinate system $(x^1\cd x^n)$
$F_f$ is given by
$$
F_f(x^1\cd x^{n-1})=(x^1\cd x^{n-1},f(x^1\cd x^{n-1})).
$$
This parameterizes a smooth hypersurface $N_f$ through $x_0$ and
moreover every smooth spacelike hypersurface through $x_0$ is uniquely
parameterized in this manner for a unique $f$ satisfying
$$
1-\sum_{i,j=1}^{n-1}g^{ij}D_ifD_jf>0.
$$
(This is exactly the condition that the image of $F_f$ is spacelike.) When
the image is spacelike set
\begin{align*}
X_i&:=\frac{\f}{\f x^i}+D_if\frac{\f}{\f x^n},\quad
W:=\biggl(1-\sum_{i,j}^{n-1}g^{ij}D_ifD_jf\biggr)^{\frac12},\\
\nor&:=\frac1W
\biggl(\frac{\f}{\f x^n}+\sum_{i,j=1}^{n-1}g^{ij}D_if\frac{\f}{\f
x^j}\biggr).
\end{align*}
Then $X_1\cd X_{n-1}$ is a basis for the tangent space to the image of $N_f$
and $\nor$ is the future pointing timelike unit normal to $N_f$.  Now a
tedious calculation 
shows that the second fundamental form $h$ of $N_f$ is given by
$$
h(X_i,X_j)=\frac{1}{W}\(D_{ij}f+\Gamma_{ij}^n-V_{ij}\)
$$
where $\Gamma^k_{ij}$ are the Christoffel symbols and $V_{ij}$ are defined by 
\begin{equation}\label{eq:Vijdef}
V_{ij} := \sum_{k=1}^{n-1}
   \(\Gamma_{ij}^kD_kf+\Gamma_{in}^kD_kfD_jf+\Gamma_{jn}^kD_kfD_if\) .
\end{equation}
Solving for the Hessian of $f$ in terms of the
second fundamental form of $N_f$ gives 
\begin{equation}\label{eq:hess2}
D_{ij}f=Wh(X_i,X_j)-\Gamma_{ij}^n+V_{ij} .
\end{equation}
The induced metric on $N_f$ has its components in the coordinate system
$x^1\cd x^{n-1}$ given by
$$
G_{ij}=g(X_i,X_j)=g_{ij} -D_ifD_jf.
$$
Let $[G^{ij}]=[G_{ij}]^{-1}$. 
Then the mean curvature of $N_f$ is
\begin{align*}
H&=\frac1{n-1}\trace_Gh=\frac1{n-1}\sum_{i,j=1}^{n-1}G^{ij}h(X_i,X_j)\\
&=\frac1{(n-1)W}\sum_{i,j=1}^{n-1}G^{ij}\(
        D_{ij}f+\Gamma_{ij}^n - V_{ij} \) ,
\end{align*}

where $x:=(x^1\cd x^{n-1})$,
$$
[ G^{ij}(x,f,Df) ] =[g_{ij}(x,f)-D_if(x)D_jf(x)]^{-1} 
$$
and
$V_{ij}(x,f,Df)$ is given by (\ref{eq:Vijdef}).
Now we can write 
$$
H  =\sum_{i,j=1}^{n-1}a^{ij}(x,f,Df)D_{ij}f +b(x,f,Df).
$$
where $a^{ij}$ and $b$ are given by 
\begin{align*}
a^{ij}(x,f,Df)&:=\frac1{(n-1)W}G^{ij}(x,f,Df),\\
b(x,f,Df)&:= \frac1{(n-1)W}\sum_{i,j=1}^nG^{ij}
\(\Gamma_{ij}^n - V_{ij} \)
\end{align*}
Therefore if ${\mathcal H}[f]$ is the mean curvature of $N_f$ then the
operator $f\mapsto {\mathcal H}[f]$ is quasi-linear.

\subsection{Proof of Theorem \ref{geo:max}}\label{sec:proofgeomax}
Given an indexed set of functions $\{f_\alpha\}$ with domains
$\Omega_{\alpha}$ so that the hypersurfaces
$N_{f_\alpha}$ are all spacelike, consider the fields of future
pointing unit timelike normals to these hypersurfaces:
$$
\nor_\alpha
  =\frac{\f}{\f x^n}
        +\frac1{W_\alpha}\sum_{i,j=1}^{n-1}g^{ij}D_if_\alpha\frac{\f}{\f x^j},
\quad\text{where}\quad
W_\alpha
   :=\biggl(1-\sum_{i,j}^{n-1}g^{ij}D_if_\alpha D_jf_\alpha \biggr)^{\frac12}
$$
The following will be used to guarantee the mean curvature operator on
the support functions to $N_0$ and $N_1$ is uniformly elliptic in the
sense of Definition~\ref{uniform}.  The proof is left to the reader.

\begin{lemma}\label{cpt:vectors}
Let $\cup_\alpha\Omega_\alpha \subset K$ where $K$ is compact.  Then
there is a compact subset $\widehat{K}$ of the tangent bundle $T(M)$
that contains the set $\cup_\alpha\{\nor_\alpha(x):x\in
\Omega_\alpha \}$ if and only if there is a 
$\rho_0>0$ so that for all $\alpha$ the lower bound
$W_\alpha(x)\ge\rho_0$ holds for $x\in \Omega_\alpha$.  Moreover if
this lower bound holds and $0<\rho<\rho_0$, there is a bound
$|f_\alpha|< B$, and if\/ $U=U_{\rho,B,K}\subset
\Re^{n-1}\times\Re\times \Re^{n-1}$ is defined by
\begin{multline}\label{def:U}
U=U_{\rho,B,K}:=\biggl\{(x^1\cd x^{n-1},r,p_1\cd p_{n-1})=(x,r,p):\\
  x\in K, |r| < B, \sum_{i,j=1}^{n-1}g^{ij}(x,r)p_ip_j < 1-\rho^2\biggr\},
\end{multline}
then for any $x\in K$ the fiber $U_x=\{((r,p):(x,r,p)\in U\}$ is
convex and for all $\alpha$, the functions 
$f_\alpha$ are $U$--admissible over $\Omega_{\alpha}$.
Finally, the mean curvature operator $\mathcal H$ is uniformly elliptic on
$U$.~\qed
\end{lemma}


Let $N_0$ and $N_1$ be as in the statement of the geometric form of the
maximum principle, Theorem~\ref{geo:max}.
Choose a coordinate system
$(x^1\cd x^n)$ centered at $q_0\in N_0\cap N_1$ that puts the metric of
$(M,g)$ in the
form~(\ref{eq:metric}).  By choosing the coordinate neighborhood
sufficiently small
we insure that, within this neighborhood, $N_0$ and $N_1$ are acausal and
edgeless, and $N_0$ is in the
causal future of $N_1$.

Then there is an open connected set $\Omega\subseteq \Re^{n-1}$ and continuous
functions $u_0, u_1:\Omega\to \Re$ so that $N_0=N_{u_0}$ and
$N_1=N_{u_1}$ as above (that is $N_i=\{(x^1\cd x^{n-1},u_i(x^1\cd
x^{n-1})): (x^1\cd x^{n-1})\in \Omega\}$).  Then $u_1\le u_0$ as $N_0$
is in the causal future of $N_1$.  Define $U:=U_{\rho_2,B,\Omega}$ by
equation~(\ref{def:U}), where $\rho_2$ is to be chosen shortly.  Then,
using that $N_1$ has mean curvature $\ge 0$ in the sense of support
hypersurfaces with one--sided Hessian bounds, the discussion above
(cf.~Lemma~\ref{cpt:vectors} and by possibly making $\Omega$ a little
smaller) there is a choice of $\rho_2\in (0,1)$ and $B>0$ so that
every $C^2$ past support hypersurface $P$ to $N_1$ has its unit
normal at the point of contact with $N_1$ in the set
$U=U_{\rho_2,B,\Omega}$ and the closure of $U$ in $T(M)$ is compact.
Moreover, from the calculations of \S\ref{reduce}, we see that any
$C^2$ past support hypersurface $P_{q,\e}$ will be of the form
$P_{q,\e}=N_{\psi_{q,\e}}$ where $\psi_{q,\e}$ is a $C^2$
$U$--admissible lower support function to $u_1$ at $q$ that satisfies
${\mathcal H}[\psi_{q,\e}](q)\ge H_0-\e$ ($\mathcal H$ is the mean curvature
operator) and which has a lower bound on its Hessian. Thus in $\Omega$
we see ${\mathcal H}[u_1]\ge H_0$ in the sense of support functions with a
one--sided Hessian bound.  Thus ${\mathcal H}$ and $u_1$ satisfy the
hypothesis of the analytic maximum principle.

Likewise if $S$ is a $C^2$ future support hypersurface for $N_0$
then there is a $C^2$ upper support function $\phi_0$ for $u_0$
such that locally $S = N_{\phi_0}$.
The problem, as noted in
Remark~\ref{geo-prob}, is that $\phi_0$ need not be $U$--admissible.  Now
use the notation of the standard setup (with $n$ replaced by
$n-1$) except that we do not know that the upper support function
$\phi_0$ to $u_0$ at $x_*$ is $U$--admissible.  As the set $U$ is open
and the function $u_0$ is continuous there is an $r_2>0$ and
$\delta_2>0$ so that
\begin{equation}\label{U-open}
\|x_*-x_1\|\le r_2,\ \text{and}\,\ \|p-D\phi_1(x_*)\|\le \delta_2
\quad\text{implies}\quad (x_*,u_0(x_*),p)\in U
\end{equation}
We now make a few modifications in \S\ref{proof}, where the proof of the
analytic maximum principle Theorem \ref{thm:max} was given.  
Change the definition of $r_1$ to
$$
r_1=\min\left\{r_0,r_2,\
\frac1{4\(n^2\ecst(\hcst+1)\)\(\ol{\alpha} r_0^{-(\ol{\alpha}+2)}\)}\right\}
$$
where $\ol{\alpha}$ is given by
$$
\ol{\alpha}:=-2+\ecst\(1+(n-1)\ecst+(n-1)^3\ecst(\hcst+1)\).
$$
(This is equation~(\ref{eq:alpha-bar}) with $n$ replaced by $n-1$.)
Then $\ol{\alpha}$ only depends on $\ecst$, $\scst$, $r_0$, and $n-1$.

If $f=\phi_1-\phi_0+\delta w$ has a local maximum at $x_*$ then by
Lemma~\ref{grad-diff}
$$
\|D\phi_1(x_*)-D\phi_0(x_*)\|\le \delta \alpha r_0^{-(\alpha+2)}.
$$
Therefore define
$$
\delta_3=\delta_2\bigl(\ol{\alpha}r_0^{-(\ol{\alpha}+2)}\bigr)^{-1}
$$
so that $\delta_3\(\ol{\alpha}r_0^{-(\ol{\alpha}+2)}\)=\delta_2$.
Now in the proof of the analytic maximum principle change the  definition of
$\delta$ to
$$
\delta=\min\{\delta_1,\delta_3,\ol{\delta}(\ol{\alpha})\}
$$
where $\ol{\delta}(\alpha)$ is still defined by equation~(\ref{r-delta})
and $\delta_1$ has the same definition as in \S\ref{proof}.  Then using
these definitions of $r_1$ and $\delta$ in the proof along with the fact
that the bound
just given on $\|D\phi_1(x_*)-D\phi_0(x_*)\|$,
implication~(\ref{U-open}) and  the definition of $\delta_3$ imply that, as
$f=\phi_1-\phi_0+\delta w$ has a local maximum at $x_*$, $\phi_0$
will be $U$--admissible. 
The rest of the proof proceeds exactly as
before to show that $u_0=u_1$ in $\Omega$.
This implies that $N_0$ and $N_1$ agree in a neighborhood ${\mathcal O}$
of $q_0$.


Finally note that as the metric $g$ is smooth the functions $a^{ij}$
and $b$ in the definition of the mean curvature operator $\mathcal H$ are
$C^\infty$.  Thus the regularity part of Theorem~\ref{thm:max} implies that
$N_1\cap {\mathcal O} =N_2 \cap {\mathcal O}$ is a smooth hypersurface.
This completes the proof of Theorem \ref{geo:max}.

\subsection{A Geometric Maximum Principle for Riemannian
Manifolds}\label{sec:Riemannian}


We now fix our sign conventions on the imbedding invariants of smooth
hypersurfaces in a Riemannian manifold $(M,g)$.  It will be convenient
to assume that our hypersurfaces are the boundaries of open sets.  As
this is always true locally it is not a restriction. As in the
Lorentzian case we denote the metric connection by $\nabla$.  Let
$D\subset M$ be a connected open set and let $N\subset\f D$ be part of
the boundary that is a $C^2$ hypersurface (we do not want to assume
that all of $\f D$ is smooth).  Let $\nor$ be the outward pointing unit
normal along $N$.  Then the second fundamental form of $N$ is the
symmetric bilinear form defined on the tangent spaces to $N$ by $
h^N(X,Y)=\la \nabla_XY,\nor\ra.  $ The mean curvature of $N$ is then $
H^N:=\frac{1}{n-1}\trace_{g\big|_N}h^N
=\frac{1}{n-1}\sum_{i=1}^{n-1}h^N(e_i,e_i) $ where $e_1,\dots,e_{n-1}$
is a local orthonormal frame for $T(N)$.  This is the sign convention
so that for the boundary $S^{n-1}$ of the unit ball $B^n$ in $\Re^n$
the second fundamental form $h^N=-g\big|_{S^n}$ is negative definite
and the mean curvature is $H^{S^{n-1}}=-1$.

\begin{definition}\label{mean-hyper-def}
Let $U$ be an open set in the Riemannian manifold $(M,g)$. Then 
\begin{enumerate}
\item   $\f U$ has mean curvature  $\ge H_0$ {\bf in the sense of contact
hypersurfaces\/} iff for all $q\in \f U$ and $\e>0$ there there is an
open set $D$ of $M$ with $\ol{D}\subseteq \ol{U}$, $q\in \f D$, and
the part of $\f D$ near $q$ is a $C^2$ hypersurface of $M$ and at the
point $q$, $H^{\f D}_q\ge H_0-\e$.
\item $\f U$ has mean curvature $\ge H_0$ {\bf in the sense of contact
hypersurfaces with a one sided Hessian bound\/} iff for all compact
$K\subset \f U$ there is constant $C_K\ge 0$ so that for all $q\in K$
and $\e>0$ there is a an open set $D$ of $M$ with $\ol{D}\subseteq
\ol{U}$, $q\in \f D$, the part of $\f D$ near $q$ is a $C^2$
hypersurface of $M$ and at the point $q$, $H^{\f D}_q\ge H_0-\e$ and
also $h^{\f D}_q\ge -C_K g\big|_{\f D}$.~\qed
\end{enumerate}
\end{definition}

The proof of the Lorentzian version of the geometric maximum principle
can easily be adapted to prove the following Riemannian version.

\begin{thm}[Geometric Maximum Principle for Riemannian Manifolds]
\label{Riem-max}
Let $(M,g)$ be a Riemannian manifold, $U_0,U_1\subset M$ open sets, and let
$H_0$ be a constant.  Assume that
\begin{enumerate}
\item $U_0\cap U_1=\emptyset$,
\item $\f U_0$ has mean curvature $\ge -H_0$ in the sense of contact
hypersurfaces,
\item $\f U_1$ has mean curvature $\ge H_0$ in the sense of contact
hypersurfaces with a one sided Hessian bound, and
\item there is a point $p\in \overline{U}_0\cap \overline{U}_1$ and a 
neighborhood $\mathcal N$ of $p$ that has coordinates $x^1,\dots,
x^n$ centered at $p$ so that for some $r>0$ the
image of these coordinates is the box $\{(x^1,\dots, x^n): |x^i|<0\}$
and there are Lipschitz continuous functions
$u_0,u_1:\{(x^1,\dots,x^{n-1}): |x^i|< r\}\to (-r,r)$ so that
$U_0\cap\mathcal{N}$ and $U_1\cap \mathcal{N}$ are given by
\begin{align*}
U_0\cap \mathcal{N}
	&=\{(x^1,\dots,x^n) : x^n >u_0(x^1,\dots,x^{n-1})\},\\ 
U_1\cap \mathcal{N}
	&=\{(x^1,\dots,x^n)  : x^n <u_1(x^1 ,\dots,x^{n-1})\}.
\end{align*}
(This implies $u_1\le u_0$ and $u_1(0,\dots,0)=u_0(0,\dots,0)$).
\end{enumerate}
Then $u_0\equiv u_1$ and $u_0$ is a smooth function.  Therefore $\f
U_0\cap \mathcal{N}=\f U_1\cap\mathcal{N}$ is a smooth embedded hypersurface
with constant mean curvature $H_0$ (with respect to the 
outward normal to $U_1$).~\qed
\end{thm}

In proving this, note that the Lipschitz conditions on $u_0$ and $u_1$
makes the mean curvature operator uniformly elliptic in the sense of
Definition~\ref{uniform} when applied to upper and lower support
functions to $u_0$ and $u_1$.

\begin{remark} \label{rk:Riem}
This is a special case of a more general result which we give in
another paper~\cite{Howard:inner-sphere}.  There the conditions on the
boundaries $\f U_i$ are relaxed to assuming only that at each point
$p$ there is an open ball $B\subset U_i$ that has $p$ in its closure
and that the radii of these balls is locally bounded from below (this
allows sets where the boundaries do not have to be topological
manifolds).  Most of the work in proving the more general maximum
principle involves proving a structure theorem for the boundaries of
open sets that satisfy a ``locally uniform inner sphere condition'' as
above, but once this is done a main step in the proof is exactly the
version of the maximum principle given in Theorem~\ref{Riem-max}.~\qed
\end{remark}

\section{Applications to Warped Product Splitting Theorems}
\label{sec:appl:warp}

\subsection{Statement of Results}
\label{sec:warpstate}
The Lorentzian geometric maximum principle, Theorem~\ref{geo:max}, provides an
especially natural and conceptually transparent proof of the
Lorentzian splitting theorem.  Recall, the Lorentzian splitting
theorem asserts that if a globally hyperbolic or timelike geodesically
complete spacetime $M$ obeys the strong energy condition,
$\operatorname{Ric}(T,T) \ge 0$ for all timelike vectors $T$, and
contains a complete timelike line $\gamma : (-\infty,\infty) \to M$
then $M$ splits isometrically into a Lorentzian product.  The proof
makes use of Lorentzian Busemann functions (see e.g.,
\cite{Beem-Ehrlich:Lorentz2} for a nice introduction).  Let $b^+$ be
the Busemann function associated to the ray $\gamma|_{[0,\infty)}$,
and let $b^-$ be the Busemann function associated to the ray
$-\gamma|_{[0,\infty)}$.  The regularity theory of Lorentzian Busemann
functions (cf.~\cite{Galloway-Horta} and
\cite[\S14.1-\S14.3]{Beem-Ehrlich:Lorentz2})
guarantees that $b^{\pm}$
are continuous on a neighborhood ${\mathcal O}$ of $\gamma(0)$ and that
$N^{\pm}= \{b^{\pm}= 0\}$ are $C^0$ spacelike hypersurfaces in the sense of
Definition \ref{def:spacelike} in ${\mathcal O}$.
$N^{\pm}$ both pass through $\gamma(0)$, and, by the reverse triangle
inequality,
$N^-$ is locally to the future of $N^+$ near $\gamma(0)$.  The curvature
assumption implies that $N^-$ has mean curvature $\le 0$ in the sense of
support
hypersurfaces and $N^+$ has mean curvature $\ge 0$ in the sense of support
hypersurfaces with one--sided Hessian bounds.  Thus, Theorem~\ref{geo:max}
implies that $N^+$ and $N^-$ agree and are smooth near $\gamma(0)$.  It follows
that $b^+ = b^- = 0$ along a smooth maximal (i.e., mean curvature zero)
spacelike hypersurface $N$.  The original proof of this given in
\cite{galloway:split} is less direct and less elementary, as it makes use
of a deep
existence result for maximal hypersurfaces due to
Bartnik~\cite{Bartnik:variational}.
It is then straight-forward to show that the normal exponential map along
$N$ gives a splitting of a tubular neighborhood of $\gamma$.  This local
splitting
can then be globalized.

The following theorem extends the Lorentzian splitting theorem to
spacetimes which
satisfy the strong energy condition with positive cosmological constant.

\begin{thm}[Warped product splitting theorem]\label{warpspli}
Let $(M,g)$ be a connected globally hyperbolic spacetime which satisfies
$\Ric(T,T) \geq n-1$ for any timelike unit vector $T$. Assume there is a
timelike  arclength parameterized geodesic segment
$\gamma : (-\pi/2 , \pi/2 ) \to M $ which maximizes the distance
between any two of its points. Then there is a complete Riemannian
manifold $(N,g_N)$ such that $(M ,g)$ is
isometric to a warped product $(-\pi/2 , \pi/2) \times N$ with metric
\begin{equation}
\label{warppromet}
-dt^2 + \cos(t)^2 g_N.
\end{equation}
\end{thm}

\begin{remark}
A calculation shows that a spacetime with metric of the form
(\ref{warppromet})
satisfies the assumptions of Theorem~\ref{warpspli} if the sectional
curvature of $(M,g)$ is bounded from below by $-1$.~\qed
\end{remark}

The following
notion which is weaker than global hyperbolicity was used by Harris in
\cite{Harris:max}.
\begin{definition}\label{piglohyp}
Let $(M,g)$ be a Lorentzian manifold and let $q > 0$.
Then $(M,g)$ is globally hyperbolic of
order $q$ if and only if  $M$ is strongly causal and
$x \ll y$, $d(x,y) < \pi/q$
implies that $C(x,y)$ is compact, where
$C(x,y)$ is the set of causal curves connecting $x$ and $y$.
~\qed
\end{definition}

To state our next result we need some notation.  Let
$(\hyps^{n-1},g_\hyps)$ be the standard Riemannian hyperbolic space with
sectional curvature $\equiv-1$ and let $(\els^n_1(-1),g_{\els})$ be
the Lorentzian manifold $\els^n_1(-1):=\hyps^{n-1}\times(\pi/2,\pi/2)$
with the metric $g_{\els}:= \cos^2(t)g_\hyps-dt^2$.  This is a
Lorentzian manifold (the subscript in $\els^n_1(-1)$ is to indicate
the metric has index one, i.e. it is Lorentzian) which has
sectional curvature~$\equiv-1$. Recall that the unique simply
connected geodesically complete Lorentzian manifold with sectional
curvature $\equiv -1$, which we will denote by $\Re^n_1(-1)$, is the
universal anti-de~Sitter space
(cf.~\cite[p.~183]{Beem-Ehrlich:Lorentz2}).
While $(\els^n_1(-1),g_{\els})$ has sectional
curvature $\equiv-1$ it is not geodesically complete and so is not
isometric to $\Re^n_1(-1)$ but is isometric to an open subset of
$\Re^n_1(-1)$.  If $x\in \hyps^{n-1}$ then the curve
$\gamma_x:(\pi/2,\pi/2)\to \els^n_1(-1)$ given by $\gamma_x(t)=(x,t)$
is a timelike unit speed geodesic and moreover
$(\els^n_1(-1),g_{\els})$ is globally hyperbolic and if $x$ is any
point of $\els^n_1(-1)$ then there are points $p$ and $q$ on
$\gamma_x$ so that $p\ll x\ll q$.  Using these facts it is not
hard to show that given any unit speed timelike geodesic
$\gamma:(\pi/2,\pi/2)\to \Re^n_1(-1)$ of length~$\pi$, the set
$\{x\in \Re^n_1(-1): \text{there are $p$ and $q$ on $\gamma$ so that
$p\ll x\ll q$}\}$ is isometric to $(\els^n_1(-1),g_{\els})$.

\begin{cor}[Lorentzian maximal diameter theorem]\label{Lor-max-diam}
Let $(M ,g)$ be a connected
Lorentzian manifold which is globally hyperbolic of order 1 and assume
that $\Ric(T,T) \geq (n-1)$ for any timelike unit vector $T$.
If $M $  contains a timelike geodesic segment
$\gamma: [-\pi/2, \pi/2] \to M $
of length $\pi$ connecting $x$ and $y$,
then 
$D = \{ z : x \ll z \ll y \}$ is
isometric to $(\els^n_1(-1),g_{\els})$.
Moreover, if $M $ contains a timelike geodesic
$\gamma: (-\infty,\infty) \to M$ such that each segment
$\gamma\big{|}_{[t,t+\pi]}$ is maximizing, then $(M ,g)$ is
isometric to the universal anti-de~Sitter space $\Re^n_1(-1)$.
\end{cor}


A very closely related result is given by
Eschenburg~\cite[Cor.~2~p.~66]{Eschenburg:max}.

Our next result makes use of the notion of the cosmological time
function (Definition~\ref{time-fcn}) which was introduced and studied
in~\cite{andersson:galloway:howard:timefcn} as a canonical choice of
a global time function for cosmological spacetimes.


\begin{definition}\label{time-fcn}
Let $(M,g)$ be a spacetime.  Then the {\bi cosmological time
function}\/ $\tau:M\to (0,\infty]$ is defined by
$$
\tau(q)=\sup\{d(p,q): p\ll q\}.
$$
where $d$ is the Lorentzian distance function.
~\qed
\end{definition}

If $(M,g)$ is spacetime of dimension $\ge4$ let $W_g$ be the Weyl
conformal tensor of $(M,g)$ written as a $(0,4)$ tensor.  Then the
squared norm of $W$ with respect to the metric $g$ is
$$
\|W_g\|^2_g = \sum g^{AA'}g^{BB'}g^{CC'}g^{DD'}W_{ABCD}W_{A'B'C'D'}.
$$
Under a conformal change of metric $\tilde{g}=\lambda^2g$ this transforms
as
\begin{equation}\label{transform}
W_{\tilde{g}}=\lambda^2W_{g},\quad
\|W_{\tilde{g}}\|^2_{\tilde{g}}=\lambda^{-4}\|W_g\|_g^2 .
\end{equation}
Because the metric $g$ is not positive definite it is possible that
$\|W_g\|_g^2=0$ at a point without $W_g$ being zero at the point and in
general $\|W_g\|_g^2$ can be negative.

\begin{thm}\label{conformal}
Let $(M,g)$ be a globally hyperbolic spacetime of dimension $n\ge4$.
Assume
\begin{enumerate}
\item $\ric(T,T) \ge (n-1)$ on all timelike unit vectors $T$.
\item There is a geodesic $\gamma:(-\pi/2,\pi/2)\to M$ that maximizes
the distance between any two of its points.
\item The Weyl conformal tensor $W_g$ of $g$ and the cosmological time
function $\tau:M\to (0,\infty]$ satisfy
$$
\lim_{\tau(q)\searrow 0}\tau(q)^4\|W_g\|_g^2=0.
$$
\end{enumerate}
Then there is a complete Riemannian manifold $(N,g_N)$ of constant
sectional curvature such that $(M ,g)$ is isometric to a warped
product $(-\pi/2 , \pi/2) \times N$ with metric
$$
-dt^2 + \cos(t)^2 g_N.
$$
\end{thm}

\begin{remark} 
This result is loosely related to the Weyl curvature
hypothesis of R. Penrose, see \cite{tod:conformal}.
Theorem \ref{conformal}  can be restated as saying that if
$\ric(T,T)\ge n-1$ on timelike unit vectors $T$, there is a line of
length $\pi$ and the Weyl conformal tensor $W_g$ has order
$o(\tau^{-2})$ (so that the squared norm $\|W_g\|_g^2$ has order
$o(\tau^{-4})$ then $(M,g)$
 is ``locally spatially isotropic''.~\qed
\end{remark}

\subsection{Proofs}
\newcommand{\bsigma}{\bar\sigma}

\subsubsection{Proof of Theorem~\ref{warpspli}} Let $(M,g)$ be an
$n$~dimensional globally hyperbolic spacetime so that $\ric(T,T)\ge(n-1)$
for all timelike unit vectors $T$.  Call a unit speed timelike
geodesic $\gamma:(-\pi/2,\pi/2)\to M$ that maximizes the distance between
an two of its points a {\bi line}\/ in $M$.  By hypothesis $M$ has at
least one line.  For any curve (or any subset of $M$) $c:(a,b)\to M$  let
$I^+(c)=\{x:\text{there is a $p$ on $c$ with $p\ll x$}\}$,
$I^-(c)=\{x:\text{there is a $q$ on $c$ with $x\ll q$}\}$ and
$I(c)=I^+(c)\cap I^-(c)$.

\begin{definition}\label{def:asy}
Let $\gamma:(-\pi/2,\pi/2)\to M$ be a line in $M$, and let 
$s \in (-\pi/2,\pi/2)$. For $p \ll \gamma(s)$, let $\alpha_s$ be a
maximal geodesic connecting $p$ and $\gamma(s)$.  If there is a
sequence $\{s_k\}_{k=1}^{\infty} \subset (-\pi/2,\pi/2)$ and a
timelike unit vector $v$ such that $s_k \nearrow \pi/2$, $p \ll
\gamma(s_k)$ and $\dot\alpha_{s_k}(0) \to v \in T(M)_p$ then the
maximal geodesic starting at $p$ in the direction $v$ is called an
{\bi asymptote}\/ to $\gamma$ at~$p$.
~\qed
\end{definition}

If $(-\pi/2,\pi/2)\to M$ is a line in $M$ then for each $r\in
(-\pi/2,\pi/2)$ define
$$
b_r(x):=d(\gamma(0),\gamma(r))-d(x,\gamma(r))=r-d(x,\gamma(r))
$$
where $d$ is the Lorentzian distance function.  The {\bi Busemann
function}\/ of $\gamma$ is then defined by
\begin{equation}\label{def-bus}
b(x):=\lim_{r\nearrow \pi/2} b_r(x).
\end{equation}
For fixed $x\in M$ the reverse triangle inequality (which is
$d(p,q)+d(q,z)\le d(p,z)$ if $p\ll q\ll z$) implies the function
$r\mapsto b_r(x)$ is monotone decreasing thus $b(x)$ exists for all
$x\in M$.  Most of the following is contained in the papers of
Eschenburg~\cite{Eschenburg:split} and Galloway~\cite{galloway:split}.

\begin{prop}\label{nice} Let $(M,g)$ be as in the statement of
Theorem~\ref{warpspli} and
$\gamma:(-\pi/2,\pi/2)\to M$ be a line in $M$.  Then there is a
neighborhood $U$ of $\gamma(0)$ (a {\bi nice neighborhood}\/) so that
\begin{enumerate}
\item The Busemann function of $\gamma$ is continuous on $U$ and satisfies
the reverse Lipschitz inequality
\begin{equation}\label{back-lip}
b(q)\ge b(p)+d(p,q)\quad\text{for}\quad p\ll q.
\end{equation}
\item For each $x\in U$ there is an asymptote $\alpha_x$ at $x$ and it
satisfies
\begin{equation}\label{ray}
b(\alpha_x(t))=b(x)+t\quad \text{for} \quad 0\le t< \pi/2-b(x).
\end{equation}
Moreover $\alpha_x$ is a ray in the sense that it is future inextendible and
maximizes the distance between any two of its points.
\item The set of initial vectors $\{\alpha_x'(0):x\in U\}$ is contained in a
compact subset of $T(M)$.
\item The zero set $N^+:=\{x\in U:b(x)=0\}$ of $b$ in $U$ is a $C^0$ 
spacelike  hypersurface (Definition \ref{def:spacelike})
in $U$ in and the mean curvature of
$B^+$ satisfies $H\ge0$ in the sense of support hypersurfaces with one-sided
Hessian bounds.
\end{enumerate}
\end{prop}

\begin{proof} Parts (1), (2), and (3) are contained in \cite{Eschenburg:split}
and \cite{galloway:split}.

It remains to prove  the bound on the mean curvature of $N^+$. Let
$x\in N^+$ and let $\eta=\alpha_x'(0)$.  Then set $r_0=\pi/4$,
$\delta=\pi/8$ and for $r\in [\pi/4,\pi/2)$ let $S_{\eta,r}$ be
defined by equation~(\ref{sphere-def}).  Then as
$\exp(t\eta)=\alpha_x(t)$ is a ray this implies, using the notation of
Proposition~\ref{sphere-bds}, that
$\gamma_\eta(t):=\exp(t\eta)=\alpha_x(t)$ maximizes distance on
$[0,r_0+\delta]$.  Also as $\alpha_x$ is a ray the segment
$\gamma_\eta\big{|}_{[0,r]}$ maximizes the distance between $x=\alpha_x(0)$ and
$\exp(r\eta)=\alpha_x(r)$ this segment contains no cut points.  This
implies the past geodesic sphere $S_{\eta,r}$ is smooth near the point
$x$, and by Proposition~\ref{sphere-bds} there will be a uniform lower
bound on the second fundamental forms of $h^{S_{\eta,r}}_x$ for $r\ge
r_0=\pi/2$ as $x$ varies over $U\cap N^+$.  By a standard comparison
theorem the mean curvature of $S_{\eta,r}$ at $x$ satisfies
$H^{S_{\eta,r}}_x\ge -\cot(r)$.  Thus for any $\e>0$ if we
choose $r\in[\pi/4,\pi/2)$ so that $\cot(r)<\e$ then $S_{\eta,r}$ will
be past support hyperplane for $N^+$ at $x$ with mean curvature at
$x>-\e$. (That $S_{\eta,r}$ is in the past $N^+$ follows from the
reverse Lipschitz inequality.) This completes the proof.
\end{proof}

We first prove a local version of Theorem~\ref{warpspli}.
Given the line $\gamma$,
let $b^+:M\to \Re$ be defined by
equation (\ref{def-bus}) above.  For $r\in (-\pi/2,0]$ let
$$
b_r^-(x):=d(\gamma(r),\gamma(0))-d(\gamma(r),x) =-r-d(\gamma(r),x)
$$
and
$$
b^-:(x)=\lim_{r\searrow-\pi/2}b_r^+(x).
$$
This is just the Busemann function for $\gamma$ if the time
orientation of $M$ is reversed and so the obvious variant of the last
proposition shows that there is a neighborhood $U$ of $\gamma(0)$ so
that if $N^\pm:=\{x\in U:b^\pm(x)=0\}$, then both $N^+$ and $N^-$ are
$C^0$ spacelike hypersurface
in $U$, the mean curvature of $N^+$ is
$\ge 0$ in the sense of support hypersurfaces with one--sided Hessian
bounds and the mean curvature of $N^-$ is $\le 0$ in the sense of
support hypersurfaces.

By the reverse triangle inequality
\begin{align}
b^+(x)+b^-(x)&=\lim_{r\nearrow\pi/2}(b^+_r(x)+b^-_{-r}(x))\nn\\
&=\lim_{r\nearrow\pi/2}(d(\gamma(0),\gamma(r))+d(\gamma(-r),\gamma(0))\nn\\
&\qquad\qquad\qquad    -d(x,\gamma(r))-d(\gamma(-r),x))\nn\\
&=\lim_{r\nearrow\pi/2}(d(\gamma(-r),\gamma(r))-d(x,\gamma(r))-d(\gamma(-r),
x))\nn\\
\label{b+-ineq}&\ge0.
\end{align}
Thus if $x\in N^-$, so that $b^-(x)=0$, then
$b^+(x)=b^+(x)+b^-(x)\ge0$.  But $\{x:b^+(x)\ge 0\}$ is in the causal
future of the set $\{x:b^+(x)=0\}$.
This implies that $N^-$ is locally to the future of
$N^+$ near $\gamma(0)$.  Therefore the hypersurfaces $N^+$ and
$N^-$ satisfy the hypothesis of the geometric maximum
principle~\ref{geo:max}.
Thus, by choosing $U$ smaller, if necessary, $N^+=N^-$ and this is a
smooth hypersurface with mean curvature $H\equiv0$.

For any point $x\in N^+$ any asymptote $\alpha_x$ to $\gamma$ must be
orthogonal to $N^+$.  This is because $N^+$ is smooth and the
support hypersurfaces we used to show $N^+$ has mean curvature $\ge
0$ in the sense of support hypersurfaces where orthogonal to
$\alpha_x$ at $x$ and as $N^+$ is smooth it must have the same tangent
planes as any of its support hypersurfaces at a point of contact.
As any asymptote to $\gamma$ from $x$ must be orthogonal to $N^+$ this
also shows that at points $x\in N^+$ that the asymptote to $\gamma$
from $x$ is unique and is given by $\alpha_x(t)=\exp(t\nor(x))$ for
$0\le t\le \pi/2$ where $\nor(x)$ is the future pointing unit normal to
$N^+$ at $x$.  From the inequality~(\ref{back-lip}) and the
identity (\ref{ray}) we see that the ray $\alpha_x$ maximizes the
distance between any of its points $\alpha_x(t)=\exp(t\nor(x))$ and $N^+$
for $t\in[0,\pi/2)$.  Thus the hypersurface $N^+$ has no focal points
along $\alpha_x$.  But $N^+$ has mean curvature $= 0$ and
$\ric(\alpha_x'(t),\alpha_x'(t))\ge n-1$.  A standard part of the
comparison theory implies that $N$ has a focal point along $\alpha_x$
at a distance $\le \pi/2$ and if equality holds then the
second fundamental form of $N^+$ vanishes at $x$ and for each $t\in
[0,\pi/2)$ the curvature tensor of $(M,g)$ satisfies
$R_{\alpha_x(t)}(X,\alpha_x'(t))\alpha_x'(t)=X$ for all vectors $X$
orthogonal to $\alpha_x'(t)$. (Cf.~\cite{Heintze-Karcher} for the
Riemannian version which is easily adapted to the Lorentzian setting.)
Note that the second fundamental form of $N^+$ vanishes at all points
and thus $N^+$ is totally geodesic. 


\begin{prop}\label{local:split}
Define a map $\Phi: N^+\times(-\pi/2,\pi/2)\to M$ by
$\Phi(x,t)=\exp(t\nor(x))$ and let  $V:=\{\exp(t\nor(x)): x\in N^+,
t\in (-\pi/2,\pi/2)\}$ be the image of $\Phi$.    Then $\Phi$ is
injective and $V$ has a warped product metric of the required type.
That is 
\begin{equation}\label{Umetric}
\Phi^*g=-dt^2+\cos(t)^2g\big|_{N^+}.
\end{equation}
Also in $V$ the Busemann functions of $\gamma$ are given by
$b^+(\Phi(x,t))= t$ and $b^-(\Phi(x,t))=-t$.  Thus $b^+$ and $b^-$ are
smooth in a neighborhood of $\gamma$ and moreover $\gamma$ is an integral
curve of the gradient field $-\nabla b^+=\partial/\partial t$ and the
level sets of $b^{+}=t$ are smooth hypersurfaces in $V$ and so the
distribution $(\nabla b^+)^\bot$ is integrable.
\end{prop}

\begin{proof} Consider the restriction $\Phi\big|_{N_+\times
[0,\pi/2)}:N_+\times [0,\pi/2)\to M$.  Then a direct calculation using
the form of the curvature along the geodesics $\alpha_x$ shows that
$\Phi^*g=-dt^2+\cos(t)^2 g_{N^+}$.  To see this map is injective
assume $\Phi(x_1,s)=\Phi(x_2,t)$ with $s,t>0$.  Then as the second
argument of $\Phi$ is the distance from $N_+$ we have $s=t$.  Now
assume that $x_1\ne x_2$ but $\Phi(x_1,s)=\Phi(x_2,s)$.  Then choose
$s_1\in (s,\pi/2)$ and let $\beta$ be the broken geodesic form $x_1$ to
$\alpha_{x_2}(s_1)$ formed by following $\alpha_{x_1}$ from $x_1$ to
$\alpha_{x_1}(s)=\alpha_{x_2}(s)$ and then following $\alpha_{x_2}$
from $\alpha_{x_1}(s)=\alpha_{x_2}(s)$ to $\alpha_{x_2}(s_1)$.  The
length of $\beta$ is $s_1$, but the corner of $\beta$ can be smoothed
to find a curve $\beta_1$ near $\beta$ that has length greater than
$s_1$.  But this contradicts that $\alpha_{x_2}$ realizes the distance
between $N_+$ and its point $\alpha_{x_2}(s_1)$ and so
$\Phi\big|_{N_+\times [0,\pi/2)}$ is injective.  
Therefore the set $V^+:=\{\exp(t\nor(x)): x\in N^+, 0\le t<\pi/2\}$
has a warped product metric of the required type.  Similar
considerations working with $N^-$, $b^-$, and the past pointing unit
normal to $N^-$ shows that the set $V^-:=\{\exp(t\nor(x)): x\in N^+,
-\pi/2<t\le 0\}$ also has a warped product metric of the required
type.  But as $N^+=N^-$ the two sets $V^+$ and $V^-$ piece together to
give that $V:=\{\exp(t\nor(x)): x\in N^+, -\pi/2<t\le \pi/2\}$ has a
warped product product metric as claimed.  The facts about the
Busemann functions are now clear.~\end{proof}

%

All that remains is to show that local warped product splitting implies
global warped product splitting.

\begin{definition}
A {\bi strip}\/ is a totally geodesic immersion $f$ of
$$
\left ( ( -\pi/2, \pi/2) \times I , - dt^2 + \cos^2(t) ds^2 \right )
$$
into $M$ for some interval $I$ so that
$f \big{|}_{(-\pi/2,\pi/2)\times \{s\}}$ is a timelike line for each $s \in
I$.
We will denote by $S$ the space $(-\pi/2,\pi/2) \times \Re$ and by $g_S$ the
metric $-dt^2 + \cos^2(t) ds^2$.
~\qed
\end{definition}

We will make use of the following elementary facts.
\begin{lemma}
\label{LemmaA}
Consider $(S,g_S)$ as just defined and let $d_q(x) = d(q,x)$
for $q,x$ causally related. Then for $0<a<\pi$
$$
\{ d_{(\pi/2,0)} = a \} = \{ \pi/2 -a \} \times \Re.
$$
Similarly, let $\gamma_1(t) = (t,s_1)$ and $\gamma_2(t) = (t,s_2)$ be two
lines. Then
$$
\lim_{t \to \pi/2} d(\gamma_1(t),\gamma_2(\tau)) = \pi/2 - \tau .\qed
$$
\end{lemma}

The following is analogous to \cite[Prop. 7.1]{Eschenburg:split},
see also  the remarks on p.~384 of  \cite{galloway:split}.
\begin{prop}\label{prop:strip}
Let $\gamma: (-\pi/2,\pi/2) \to M$ be a timelike
line and let $\sigma: [0,\ell] \to M$ be a geodesic
with $\sigma(0) = \gamma(0)$. Then there is a strip containing $\gamma$ and
$\sigma$.
\end{prop}
\begin{proof}
Suppose there is a strip $f: (-\pi/2,\pi/2 ) \times [0,a) \to M$ containing
$\gamma$ and $\sigma\big{|}_{[0,v)}$ so that for all $t$,
$\gamma(t) = f(t,0)$, where $v$ is the defined as
$$
v = \sup \{ u : \sigma\Big{|}_{[0,u)} \subset f( (-\pi/2,\pi/2 ) \times [0,a) )
\} .
$$
Now consider the geodesic $\lambda$ in the space $U := (-\pi/2 , \pi/2 )
\times \Re$ with the metric~(\ref{Umetric}) with
$f_{*} \dot \lambda(0) = \dot \sigma (0)$. Define $X_0 \in T(N)_{(0,0)}$ by
parallel translating $\partial/\partial t$ from $\lambda(v)$ to $\lambda(0)$
along $\lambda$. Then define a vector $X \in T(M)_{\sigma(v)}$ by parallel
translating $f_* X_0$ along $\sigma$ from $\sigma(0)$ to $\sigma(v)$. By
construction, $X$ is a timelike unit vector.

Let $(t_0, s_0)$ be the coordinates of $\lambda(v)$ and let
$\gamma_v: (a,b) \to M$ be
the inextendible arclength parameterized timelike geodesic such that
$\dot \gamma_v(t_0) = X$.

We claim that $(a,b) = (-\pi/2 , \pi/2)$. Assume that this is not the case,
for example $b < \pi/2$. Then, $\gamma_v(t)$ for $t \in
[t_0,b)$ are limits of points of the form
$$
f(t,s), \ (t,s) \in [t_0,b)\times [0,s_0) .
$$
Note that $( ( -\pi/2, \pi/2 ) \times I, - dt^2 + \cos^2(t) ds^2)$ is
conformal to $(\Re\times I, - {dt'}^2 + ds^2)$ 
(by the change of coordinates $t' = \ln(\sec(t) + \tan(t))$ for 
any interval $I$.
Since the causal structure is invariant under conformal changes,
we see that we may choose $\tau < \pi/2$ sufficiently large, so that
$f([0,b)\times [0,s_0)) \subset J^{-}(\gamma(\tau))$. Therefore we find that
$$
\gamma_v\Big{|}_{[t_0,b)} \subset J^{+}(\gamma_v(t_0)) \cap J^{-}(\gamma(\tau))
.
$$
But by assumption, $\gamma_v$ is future inextendible so this contradicts the
global hyperbolicity of $M$.

Therefore $b=\pi/2$ and similarly we find that $a = -\pi/2$. Further,
$\gamma_v$ is a limit of lines of the form $f(\cdot,s_n)$ for some
sequence $s_n$ converging to $s_0$ and thus $\gamma_v$ is a line. Now
applying the local splitting (Proposition~\ref{local:split}) at
$\gamma_v$ we get a contradiction to the choice of $v$. This shows
that $v=\ell$ and thus that $\sigma$ is contained in the strip.
\end{proof}

Following \cite{Eschenburg:split} we say that two lines of length $\pi$
are {\bi strongly parallel}\/ if they bound a strip and {\bi parallel}\/
 if there is a
finite sequence of lines
$\gamma_1 = \sigma_1, \sigma_2 , \dots , \sigma_k = \gamma_2$
such that $\sigma_j,\sigma_{j+1}$ are strongly parallel.
The following is analogous to \cite[Lemma 7.2]{Eschenburg:split}.
\begin{lemma}
If $\gamma_1$ and $\gamma_2$ are parallel  lines, then
$I(\gamma_1) = I(\gamma_2)$ and the Busemann functions $b_1^{\pm}$ and
$b_2^{\pm} $ of $\gamma_1$ and $\gamma_2$ agree. Thus for a line $\gamma_1$
and a point $x\in M$ there is at most one  
line $\gamma_2$ through $x$ and
parallel to~$\gamma_1$.
\end{lemma}
\begin{proof}
It is clear that it is sufficient to consider only the case of strongly
parallel lines.
It follows from the causal structure of $(S,g_S)$
that $\gamma_1 \subset I(\gamma_2)$ and similarly that $\gamma_2 \subset
I(\gamma_2)$. Therefore $I(\gamma_1) = I(\gamma_2)$.

Further, as $t \to \pi/2$, we have
$$
\limsup_{t \to \pi/2 } d ( f( 0,s) , \gamma_i(t)) \geq \pi/2
$$
by Lemma~\ref{LemmaA}. Thus, for $i=1,2$,
$b_i^{+} ( f(s,0)) \leq 0$ and $b_i^{-}(f(s,0)) \leq 0$. Therefore by
(\ref{b+-ineq}) we get $b_i^{\pm} (f(0,s)) = 0$.

This shows using (\ref{back-lip})
that $b_1^+(\gamma_2(s)) \geq s$. On the other hand using Lemma~\ref{LemmaA}
and the reverse triangle inequality we have $b_1^{+} (\gamma_2(s)) \leq s$.
Therefore we find that
$\gamma_1$ is a coray of $\gamma_2$. Thus as in \cite[(5)]{Eschenburg:split}
we have  $b_1^{\pm} \geq b_2^{\pm} \geq b_1^{\pm}$.
Thus $b_1^{\pm}=b_2^{\pm}$ as claimed.

Let $\gamma_2$ and $\gamma_3$ be lines of $(M,g)$ through $x$ and parallel
to $\gamma_1$ and let $b_2^+$ and $b_3^+$ be the Busemann functions of
$\gamma_2$ and $\gamma_3$.  Then we have just shown that $b_2^+=b_3^+$.  By
the local version, Proposition~\ref{local:split},
of the splitting we see that $b_2^+$ and $b_3^+$ are
smooth functions and $\gamma_2$ and $\gamma_3$ are both integral curves for
the vector field $-\nabla b_2$.  Thus the uniqueness of integral curves
implies $\gamma_2=\gamma_3$.   This completes the proof.
\end{proof}

Let $\gamma$ be a fixed line and let $P_{\gamma} \subset M$ be the set of
points which lie on a parallel line. The following Lemma is proved
exactly as in \cite{Eschenburg:split}.
\begin{lemma}[\protect{\cite[Lemma 7.3]{Eschenburg:split}}]
\label{LemmaB}
$P_{\gamma}$ is a connected
component of $M$.~\qed
\end{lemma}

We now complete the proof of the theorem along the lines of
Eschenburg~\cite{Eschenburg:split}.
The uniqueness of the line $\gamma_q$ through each $q \in M$ gives a
timelike vector field $V$ on $M$. Let $V^{\perp}$ be the distribution
orthogonal to $V$. It follows from 
the local version Proposition~\ref{local:split}
this vector field is smooth and $V^{\perp}$ is integrable and 
the integral leaves can be represented by the level sets of the time
function $t$.

Let $N$ be the maximal leaf through $\gamma(0)$ and let $g$ be the induced
metric on $N$. Then the map
$$
j: (-\pi/2, \pi/2) \times N \to M; \qquad j(t,q) \to \gamma_q(t)
$$
is an isometry with respect to the metric
\begin{equation}
\label{globmetric}
-dt^2 + \cos^2(t) g .
\end{equation}
The local part of this follows from the local result
Proposition~\ref{local:split} and
the global statement follows from Lemma~\ref{LemmaB}. Thus in particular
$(-\pi/2,\pi/2) \times M$ with the metric (\ref{globmetric}) is globally
hyperbolic. Therefore, by \cite[Theorem~3.66~p.~103]{Beem-Ehrlich:Lorentz2}
$(H,h)$ is a complete Riemannian manifold.
This completes the proof of Theorem~\ref{warpspli}~\qed

\subsubsection{Proof of Corollary~\ref{Lor-max-diam}}
To prove the first statement apply Theorem~\ref{warpspli} to the interior of
$D$. This implies there is a complete Riemannian manifold $(N,g_N)$ so
that the restriction of the metric  to $D$ is a warped product
$g=-dt^2+\cos(t)^2g$ on $N\times(-\pi/2,\pi/2)$.
We now argue that because the metric $g$
is smooth at $x$ that the manifold $(N,g_N)$ must be the hyperbolic
space $(\Re^{n-1}_0(-1),g)$.  It is enough to show $(N,g_N)$ has constant
sectional curvature~$-1$.  Let $N_t:= N\times\{t\}$.  Then the
induced metric on $N_t$ is $\cos^{2}(t)g$.  Thus
the intrinsic sectional curvatures of $N_t$ and $N_0$ on corresponding
$2$~planes are related by $K_t=(1/\cos^2(t))K_0g$.  The distance of $N_t$
to $x=\gamma(-\pi/2)$ is $\pi/2+t$.  Thus
$N_t:=\{z:x\ll z,\ d(x,z)=\pi/2+t\}$.  But in the limit as
$t\searrow-\pi/2$ the set $N_t$ is then the ``geodesic sphere''
of radius $r=\pi/2+t$.  As the radius goes to $0$
(i.e. as $t\to -\pi/2$) the Gauss equation implies the sectional
curvatures behave like  $K_t=-1/r^2+{O}(1/r)
=-1/(\pi/2+t)^2+{O}(1/(\pi/2+t))$. Therefore
$$
K_t=\frac{K_0}{\cos^2(t)}=\frac{K_0}{(\pi/2+t)^2}+{ O}
\(\frac{-1}{(\pi/2+t)}\)=\frac{-1}{(\pi/2+t)^2}+{ O}\(\frac1{(\pi/2+t)}\).
$$
This implies $K_0\equiv -1$.

The proof of the second statement can be carried out along the lines of
\cite[\S 2.2]{Harris:max}.~\qed

\subsubsection{Proof of Theorem~\ref{conformal}}
Let $(M,g)$ satisfy the conditions of Theorem~\ref{conformal}, then
$(M,g)$ also satisfies the hypothesis of Theorem~\ref{warpspli}.
Therefore $(M,g)$ splits as a warped product with metric of the
form~(\ref{warppromet}).  Writing points of $(M,g)$ as $(t,x)\in
(-\pi/2,\pi/2)\times N$ it is not hard to see that the
cosmological time function $\tau$ is given by $\tau(x,t)=t+\pi/2$.
Define a new metric conformal $\tilde{g}$ to $g$ by
$$
\tilde{g}:=\frac1{\cos(t)^2}g=-\frac{dt^2}{\cos(t)^2}+g_N=-(ds)^2+g_N
$$
where
$$
s=\int_0^t\frac{dz}{\cos(z)}=\ln(\sec(t)+\tan(t))
=\ln(\sec(\tau-\pi/1)+\tan(\tau+\pi/2)).
$$
Therefore $\tilde{g}=-ds^2+g_N$ is the product metric on
$(-\infty,\infty)\times N$.  Then using the transformation
rule~(\ref{transform}) for $\|W_{\tilde{g}}\|_{\tilde{g}}^2$ and using
that near $t=-\pi/2$, $\cos(t)=\cos(\tau+\pi/2)=\tau+O(\tau^3)$ so
that
$$
\lim_{s\to-\infty}\|W_{\tilde{g}}\|_{\tilde{g}}^2=\lim_{t\searrow
\pi/2}\cos(t)^4\|W_g\|_g^2= \lim_{\tau\searrow 0}\tau^4\|W_g\|_g^2=0.
$$
But as $\tilde{g}$ is the product metric the function
$\|W_{\tilde{g}}\|_{\tilde{g}}^2$ will be constant along any of the
lines $s\mapsto (s,x)$ in $(M,\tilde{g})$.  Therefore the limit above
implies $\|W_{\tilde{g}}\|_{\tilde{g}}^2\equiv0$ on $M$.  As $(M,g)$
is a Lorentzian manifold this in not enough to conclude
$W_{\ol{g}}\equiv 0$.  To get any conclusion at all we have to use
that $(M,\ol{g})$ is a product metric.
\begin{lemma}\label{prod-ineq}
Let $(N,g_N)$ be a Riemannian manifold of dimension at least three,
set $M:=\Re\times N$, and give $M$ the Lorentzian product metric
$g:=-dt^2+g_N$.  Let $R_{ABCD}$ be the curvature tensor of $(M,g)$ as
a $(0,4)$ tensor, $R_{AB}$ the Ricci tensor and $S$ the scalar
curvature of $(M,g)$.  Let $a$ and $b$ be real numbers
and set
\begin{multline*}
V_{ABCD}:=R_{ABCD}+a(g_{AC}R_{BD}+g_{BD}R_{AC}-g_{BC}R_{AD}-g_{AD}R_{BC}) \\
        +b(g_{AC}g_{BD}-g_{AD}g_{BC}).
\end{multline*}
Then $\|V\|_g^2\ge0$ and equality on all of $(M,g)$ implies
$(N,g_N)$ has constant sectional curvature.  In particular the Weyl
conformal tensor $W$ of $(M,g)$ is of this form so if
$\|W\|_g^2\equiv0$ then $(N,g_N)$ has constant sectional curvature.
\end{lemma}

\begin{proof}
This is a pointwise result so let $p\in M$ and choose an orthonormal
basis $e_1\cd e_n$ of $T(M)_p$ so that $e_n=\f/\f t$.  Then the vectors
$e_1\cd e_{n-1}$ are tangent to the factor $N$.  In what follows we will
use the following range of indices
$$
1\le A,B,C,D,A',B',C',D\le n,\quad 1\le i,j,k,l,i',j',k',l'\le n-1.
$$
In the basis $e_1\cd e_n$
$$
g_{AB}=g^{AB}, \quad g_{ij}=\delta_{ij},\quad g_{in}=0,\quad g_{nn}=-1.
$$
With this notation we will show
\begin{equation}\label{pos-norm}
\|V\|_g^2=\sum_{i,j,k,l}(V_{ijkl})^2 + 4\sum_{i,j}(aR_{ij}+bSg_{ij})^2\ge0.
\end{equation}
To do this split the sum for $\|V\|_g^2$ into two parts, the first
only summing over the indices $i,j,k,l,i',j',k',l'$ and the second
where at least one index in the sum is equal to $n$.
\begin{align}
\|V\|_g^2 &=\sum_{\substack{A,B,C,D\\A',B',C',D'}}
g^{AA'}g^{BB'}g^{CC'}g^{DD'}
        V_{ABCD}V_{A'B'C'D'}\nn 
\\
&=\sum_{\substack{i,j,k,l\\i',j',k',l'}}
        g^{ii'}g^{jj'}g^{kk'}g^{ll'}V_{ijkl}V_{i'j'k'l'}
+{\sum}_2\nn \\
\label{split-sum}
&=\sum_{i,j,k,l}(V_{ijkl})^2+{\sum}_2
\end{align}
We now consider the terms that occur in
${\sum}_2$.  Because the metric $g=-dt^2+g_N$ is a product metric and $n$
corresponds to the direction of $e_n=\f/\f t$ we have
$$
R_{ABCn}=0,\quad R_{An}=0.
$$
All the terms in the sum ${\sum}_2$ have at least one index equal to $n$.
By the symmetries of the curvature tensor we can consider the case where
$D=n$. Using that $(M,g)$ is a product $R_{ABCn}=0, R_{An}=0$ so
$$
V_{ABCn}=a(g_{Bn}R_{AC}-g_{An}R_{BC})+bS(g_{AC}g_{Bn}-g_{An}g_{BC}).
$$
From this we see $V_{ABnn}=0$ and likewise $V_{nnCD}=0$. Also $V_{ijkn}=0$
Thus if a term in  ${\sum}_2$ is  nonzero, exactly two of
$A,B,C,D$ are equal to $n$ and moreover the case $A=B=n$ and $C=D=n$ give
zero terms.  Consider the case $B=D=n$.  Then
$$
V_{inkn}=-aR_{ik}-bSg_{ik}.
$$
The sum over the terms in ${\sum}_2$ where $B=D=n$ is then
\begin{align*}
\sum_{\substack{i,j\\A',B',C',D'}}g^{iA'}g^{nB'}g^{kC'}g^{nD'}
        V_{ABCD}V_{A'B'C'D'}&=
        \sum_{i,j}(g^{nn})^2V_{inkn}V_{inkn}\\
&=\sum_{i,k}(aR_{ij}+bSg_{ik})^2.
\end{align*}
The calculations for the three cases $A=D=n$, $A=C=n$, and $B=C=n$ give
the same result.  As these are the only cases leading to nonzero terms in
the sum ${\sum}_2$ we have ${\sum}_2=4\sum_{i,j}(aR_{ij}+bSg_{ij})^2$.
Using this in (\ref{split-sum}) completes the verification
of~(\ref{pos-norm}).

Now assume that $\|V\|_g^2=0$ at a point.  Then
equation~(\ref{pos-norm}) implies that $aR_{ij}+bSg_{ij}=0$.
Therefore
$$
a(g_{ik}R_{jl}+g_{jl}R_{ik}-g_{jk}R_{il}-g_{il}R_{jk})=
-2b(g_{ik}g_{jl}- g_{il}g_{jk})
$$
But (\ref{pos-norm}) also implies $V_{ijkl}=0$.  Using these facts in
the definition of $V$ to solve for $R_{ijkl}$ gives
$$
R_{ijkl}=bS(g_{ik}g_{jl}- g_{il}g_{jk})
$$
for $1\le i,j,k,l\le n-1$.  But as $(M,g)=(\Re\times N,-dt^2+g_N)$ is a
product of with $(N,g_N)$ the curvature tensor $R$ of $(M,g)$ and the
curvature tensor $R^N$ of $(N,g_N)$ are related by
$R^N_{ijkl}=R_{ijkl}=bS(g_{ik}g_{jl}- g_{il}g_{jk})$.  Therefore if
$\|V\|_g^2\equiv0$ then the curvature tensor of $(N,g_N)$ is of the
form $R^N_{ijkl}=bS(g_{ik}g_{jl}- g_{il}g_{jk})$.  As the dimension of
$N$ is $\ge3$ a well known theorem of Schur
(cf.~\cite[Thm~2.2 p.~202]{Kobayashi-Nomizu:vol2})
implies $(N,g_N)$
has constant sectional curvature. This completes the proof.
\end{proof}

We now complete the proof of Theorem~\ref{conformal}.  We have shown
that under the hypothesis of the theorem the metric
$\tilde{g}=-ds^2+g_N$ has $\|W_{\tilde{g}}\|_{\tilde{g}}^2 \equiv0$.
Then the lemma implies that $(N,g_N)$ has constant sectional
curvature. The completes the proof.

{\small {\bf Acknowledgment:}  We have profited from conversations and
correspondence with Joe Fu on topics related to the maximum principle
and regularity of sub-solutions to elliptic equations. }

\providecommand{\bysame}{\leavevmode\hbox to3em{\hrulefill}\thinspace}


\begin{thebibliography}{10}

\bibitem{Alexandrov:2nd-thms}
A.~D. Aleksandrov, \emph{Some theorems on partial differential equations of
  second order}, Vestnik Leningrad. Univ. Ser. Mat. Fiz. Him. \textbf{9}
  (1954), no.~8, 3--17.

\bibitem{Alexandrov:uniqueness1}
\bysame, \emph{Uniqueness theorems for surfaces in the large. {I}}, Vestnik
  Leningrad. Univ. \textbf{11} (1956), no.~19, 5--17.

\bibitem{Alexandrov:uniqueness2}
\bysame, \emph{Uniqueness theorems for surfaces in the large. {I}{I}}, Vestnik
  Leningrad. Univ. \textbf{12} (1957), no.~7, 15--44.

\bibitem{andersson:galloway:howard:timefcn}
L.~Andersson, G.~J. Galloway, and R.~Howard, \emph{The cosmological time
  function}, preprint (1997).

\bibitem{andersson:howard:RW}
L.~Andersson and R.~Howard, \emph{Rigidity results for {R}obertson--{W}alker
  and related spacetimes}, preprint (1997).

\bibitem{Bartnik:variational}
R.~Bartnik, \emph{Regularity of variational maximal surfaces}, Acta Math.
  \textbf{161} (1988), 145--181.

\bibitem{Beem-Ehrlich:Lorentz2}
John~K. Beem, Paul~E. Ehrlich, and Kevin~L. Easley, \emph{Global {L}orentzian
  geometry}, second ed., Monographs and Textbooks in Pure and Applied
  Mathematics, vol. 202, Marcel Dekker Inc., New York, 1996.

\bibitem{Caffarelli-Cabre}
Luis~A. Caffarelli and Xavier Cabr{\'e}, \emph{Fully nonlinear elliptic
  equations}, American Mathematical Society Colloquium Publications, vol.~43,
  American Mathematical Society, Providence, RI, 1995.

\bibitem{Calabi:maximum}
E.~Calabi, \emph{An extension of {E}. {H}opf's maximum principle with an
  application to {R}iemannian geometry}, Duke Math. J. \textbf{25} (1957),
  45--56.

\bibitem{Cheng:max-diam}
S.~Y. Cheng, \emph{Eigenvalue comparison theorems and its geometric
  applications}, Math. Z. \textbf{143} (1975), 289--297.

\bibitem{Eschenburg:split}
J.-H. Eschenburg, \emph{The splitting theorem for space-times with strong
  energy condition}, J. Diff. Geom. \textbf{27} (1988), 477--491.

\bibitem{Eschenburg:max}
\bysame, \emph{Maximum principle for hypersurfaces}, Manuscripta Math.
  \textbf{64} (1989), 55--75.

\bibitem{Eschenburg-Galloway}
J.-H. Eschenburg and G.~J. Galloway, \emph{Lines in space-times}, Comm. Math.
  Phys. \textbf{148} (1992), no.~1, 209--216.

\bibitem{galloway:warsaw}
G.~J. Galloway, \emph{Some rigidity results for spatially closed spacetimes},
  Proceedings of the Minisemester on Mathematical Aspects of Theories of
  Gravitation, Banach Center, Polish Academy of Science, Warsaw, to appear.

\bibitem{galloway:split}
\bysame, \emph{The {Lorentzian} splitting theorem without the completeness
  assumption}, J. Diff. Geom. \textbf{29} (1989), 373--389.

\bibitem{Galloway-Horta}
G.~J. Galloway and A.~Horta, \emph{Regularity of {L}orentzian {B}usemann
  functions}, Trans. Amer. Math. Soc. \textbf{348} (1996), no.~5, 2063--2084.

\bibitem{Gilbarg-Trudinger}
D.~Gilbarg and N.S. Trudinger, \emph{Elliptic partial differential equations of
  second order}, Springer, Berlin, 1983, 2:nd. ed.

\bibitem{Harris:max}
S.~G. Harris, \emph{On maximal geodesic-diameter and causality in {L}orentz
  manifolds}, Math. Ann. \textbf{261} (1882), 307--313.

\bibitem{Heintze-Karcher}
E.~Heintze and H.~Karcher, \emph{A general comparison theorem with applications
  to volume extimates for submanifolds}, Ann Scient. {\'E}c. Norm. Sup.
  \textbf{11} (1978), 451--470.

\bibitem{Hopf:maximum}
E.~Hopf, \emph{Elementare {B}emerkungen {\"uber} die {L\"osungen} partieller
  {D}ifferentialgleichungen zweiter {O}rdnung vom elliptischen {T}ypus},
  Sitzungsberichte, Preussische Akademie der Wissenschaften \textbf{19} (1927),
  147--152.

\bibitem{Howard:inner-sphere}
R.~Howard, \emph{Boundary structure and the strong maximum principle for sets
  satisfying a locally uniform ball conditions in a {Riemannian} manifold},
  Preprint (1997).

\bibitem{Kobayashi-Nomizu:vol2}
S.~Kobayashi and K.~Nomizu, \emph{Foundations of differential geometry},
  vol.~2, John Wiley and Sons, 1969.

\bibitem{newman:splitting}
Richard P. A.~C. Newman, \emph{A proof of the splitting conjecture of
  {S}.-{T}.\ {Y}au}, J. Differential Geom. \textbf{31} (1990), no.~1, 163--184.

\bibitem{tod:conformal}
K.~P. Tod, \emph{Isotropic singularities}, Rend. Sem. Mat. Univ. Politec.
  Torino \textbf{50} (1992), no.~1, 69--92 (1993), Singularities in curved
  space-times.

\end{thebibliography}
\end{document}